# DYNAMIC IMPORTANCE SAMPLING FOR QUEUEING NETWORKS

By Paul Dupuis,[1] Ali Devin Sezer[2] and Hui Wang[3]

*Brown University*

Importance sampling is a technique that is commonly used to speed up Monte Carlo simulation of rare events. However, little is known regarding the design of efficient importance sampling algorithms in the context of queueing networks. The standard approach, which simulates the system using an a priori fixed change of measure suggested by large deviation analysis, has been shown to fail in even the simplest network setting (e.g., a two-node tandem network).

Exploiting connections between importance sampling, differential games, and classical subsolutions of the corresponding Isaacs equation, we show how to design and analyze simple and efficient dynamic importance sampling schemes for general classes of networks. The models used to illustrate the approach include $d$-node tandem Jackson networks and a two-node network with feedback, and the rare events studied are those of large queueing backlogs, including total population overflow and the overflow of individual buffers.

**1. Introduction.** For more than two decades, there has been a growing of interest in importance sampling, a method in which the system is simulated under a different probability distribution (i.e., change of measure), for fast simulation of rare events in queueing networks [2, 11].

The standard approach to importance sampling for queueing considers an a priori fixed and static change of measure that is suggested by large deviation analysis. This approach works well for simulating large buildups of a single/multiple server queue [1, 14]. However, there has been limited success in extending this standard heuristic to networks of queues. In even the

Received December 2005; revised January 2007.
[1]Supported in part by NSF Grants DMS-03-06070 and DMS-04-04806 and the Army Research Office DAAD19-02-1-0425 and W911NF-05-1-0289.
[2]Supported in part by NSF Grants DMS-03-06070 and DMS-04-04806.
[3]Supported in part by NSF Grant DMS-04-04806.

simplest network setting, such as a two-node tandem Jackson network, the change of measure suggested by the standard heuristic fails to be asymptotically optimal in general [10, 13] and can lead to importance sampling estimators with infinite variance [3]. This failure is in fact due to the discontinuities of the state dynamics on the boundaries of the state space. Such discontinuities are not present in the case of a single queue.

The purpose of the present paper is to present a framework under which one can systematically build efficient *dynamic* (i.e., state dependent) importance sampling schemes for simulating rare events in queueing networks. Our method heavily exploits a recently discovered connection between importance sampling and deterministic differential games [7, 8] and the role of classical subsolutions of the Isaacs equation associated with the game [9]. We demonstrate that one can construct classical subsolutions, as the mollification of the minimum of affine functions, that lead to simple and efficient importance sampling schemes.

To illustrate the main idea, we focus in much of the paper on two-node tandem Jackson queueing networks. The rare events of interest are various types of buffer overflows, including total population overflow and individual buffer overflows. Also discussed are extensions to $d$-node tandem Jackson networks (Section 4) and a two-node Jackson network with feedback (Section 5). We wish to point out that our approach can be applied to general Jackson networks and networks with more general arrival/service processes (such as Markov modulated processes), and such results will be reported elsewhere. To the best of our knowledge, the present paper is the first to provide a rigorous theoretical framework in which one can build asymptotically optimal importance sampling algorithms for rare events in networks of queues.

The paper is organized as follows. Section 2 gives a brief review of the basics of importance sampling. In Section 3, we study in detail the classical problem of total population overflow in two-node tandem Jackson networks. Extensions are discussed in Sections 4 and 5. To ease exposition, most proofs are collected in the Appendices.

**2. Basics of importance sampling.** The basic idea of importance sampling is to use a change of measure, that is, the system is simulated under a different probability distribution and the outcomes are multiplied by appropriate likelihood ratios (i.e., Radon–Nikodým derivatives) to form unbiased estimators.

We specialize to the estimation of rare event probabilities and consider a family of events $\{A_n\}$ in a probability space $(\Omega, \mathcal{F}, \mathbb{P})$ such that

$$\lim_n -\frac{1}{n} \log \mathbb{P}(A_n) = \gamma > 0.$$

In order to estimate $\mathbb{P}(A_n)$, importance sampling generates samples under a probability measure $\mathbb{Q}$ such that $\mathbb{P} \ll \mathbb{Q}$, and forms an estimator by averaging independent replications of

$$\hat{p}_n \doteq 1_{A_n} \frac{d\mathbb{P}}{d\mathbb{Q}},$$

where $d\mathbb{P}/d\mathbb{Q}$ is the Radon–Nikodým derivative. It is easy to check that this importance sampling estimator is unbiased. Its rate of convergence is determined by the variance of $\hat{p}_n$. The smaller the variance, the faster the convergence. Thanks to the unbiasedness of $\hat{p}_n$, minimizing the variance amounts to minimizing the second moment, which is

(2.1) $\qquad$ [2nd moment of $\hat{p}_n$] $= E^{\mathbb{Q}}[\hat{p}_n^2] = E^{\mathbb{P}}[\hat{p}_n].$

However, Jensen's inequality implies that

$$\limsup_n -\frac{1}{n} \log E^{\mathbb{Q}}[\hat{p}_n^2] \leq \limsup_n -\frac{2}{n} \log E^{\mathbb{Q}}[\hat{p}_n] = 2\gamma.$$

We say the importance sampling estimator is *asymptotically optimal* if

$$\liminf_n -\frac{1}{n} \log E^{\mathbb{P}}[\hat{p}_n] \geq 2\gamma.$$

Sometimes $2\gamma$ is referred to simply as the "optimal decay rate."

REMARK 2.1. The requirement of $\mathbb{P} \ll \mathbb{Q}$ is more stringent than necessary. It is sufficient that $\mathbb{P}$ be absolutely continuous with respect to $\mathbb{Q}$ on a sub-$\sigma$-algebra that contains $A_n$, in which case the likelihood ratio is defined as the Radon–Nikodým derivative of $\mathbb{P}$ and $\mathbb{Q}$ when they are restricted on this sub-$\sigma$-algebra. In this paper, the changes of measure will be applied to a sequence of i.i.d. random variables $\{Y(k)\}$, and will be restricted to the $\sigma$-algebra generated by $\{Y(k)\}$ up until the time either the buffer overflow happens or the system returns to the empty state. Note that when considered on the full $\sigma$-algebra generated by $\{Y(k)\}$, it is typical that $\mathbb{P}$ is singular with respect to $\mathbb{Q}$.

**3. Two-node tandem Jackson networks.** To illustrate the main idea of the game/subsolution approach toward importance sampling, we specialize to two-node Jackson tandem queueing networks, where the arrival process is Poisson with rate $\lambda$ and the service times are distributed exponentially with rates $\mu_1$ and $\mu_2$, respectively (see Figure 1). The system is assumed to be stable, that is, $\lambda < \min\{\mu_1, \mu_2\}$.

Suppose that the two queues share one buffer with capacity $n$, and that we are interested in the overflow probability

$p_n \doteq \mathbb{P}\{\text{network total population reaches } n \text{ before returning to } 0,$

$\hspace{4cm} \text{starting from } 0\}.$

This overflow problem was among the first to be studied in the literature on importance sampling for networks, and has served as a benchmark since then [13].

Rescaling the time variable will have no effect on $p_n$, and so without loss of generality we assume $\lambda + \mu_1 + \mu_2 = 1$. Since exchanging the order of service rates does not affect this probability [16], we further assume that $\mu_2 \leq \mu_1$. Under these conditions, we have the large deviation limit [10]

$$\lim_n -\frac{1}{n} \log p_n = \log \frac{\mu_2}{\lambda} \doteq \gamma. \tag{3.1}$$

3.1. *The standard heuristic.* Based on a heuristic application of large deviation analysis, Parekh and Walrand [13] proposed a state-independent importance sampling algorithm for estimating $p_n$, which amounted to interchanging the arrival rate and the smallest service rate. For this scheme, numerical experiments suggested good performance for a certain range of parameters [13].

A rigorous analysis of this importance sampling algorithm first appeared in [10], in which the authors showed that the algorithm is asymptotically optimal when the parameters fall into certain subset. However, it was also shown that the asymptotic optimality fails in general, such as when the two service rates are nearly equal and the arrival rate is small. A recent paper [3] extended these results and showed that this scheme can lead to estimators with infinite variance for certain values of parameters. Additional discussion on importance sampling for queueing networks can be found in [11].

3.2. *The system dynamics.* The system state can be described by the embedded discrete time Markov chain $Z = \{Z(k) = (Z_1(k), Z_2(k)): k = 0, 1, 2, \ldots\}$ defined on a probability space $(\Omega, \mathcal{F}, \mathbb{P})$, where $Z_i(k)$ is the queue length at node $i$ after the $k$th transition epoch of the network.

At times when both queues are nonempty, the increments of the Markov chain $Z$ take values in the space

$$\mathbb{V} \doteq \{v_0 = (1,0),\ v_1 = (-1,1),\ v_2 = (0,-1)\},$$

with $v_0$ corresponding to an arrival and $v_i$ to a service at node $i$ for $i = 1, 2$. On the boundary where either queue is empty, the dynamics exhibit different behaviors. Suppose that the queue at node $i$ ($i = 1, 2$) is empty. Then it is impossible for the process $Z$ to have increment $v_i$ since it will lead to negative queue size. One way to describe this discontinuity is to allow $Z$ to make fictitious jumps of size $v_i$ on the boundary, but they have to be accounted for by "pushing back" the state along the direction of constraints

$$d_i = -v_i,$$

so that the state process $Z$ stays nonnegative. See Figure 2.

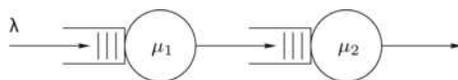

Fig. 1.  *Two-node tandem Jackson network.*

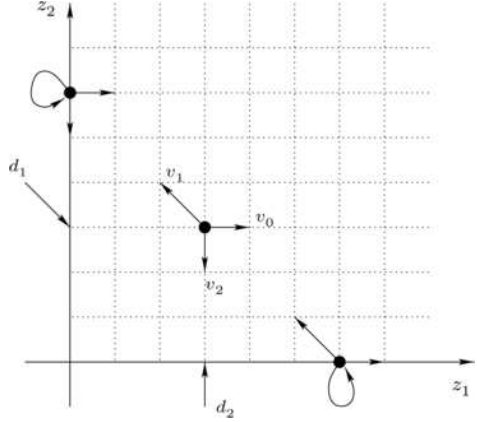

Fig. 2. *The system dynamics.*

To summarize, the evolution of the Markov chain $Z$ can be modeled by equation

$$(3.2) \qquad Z(k+1) = Z(k) + \pi[Z(k), Y(k+1)],$$

where $Y = \{Y(k) : k \geq 1\}$ is a sequence of random variables taking values in the space $\mathbb{V}$, and the mapping $\pi$ is defined for every $z = (z_1, z_2) \in \mathbb{R}_+^2$ and $y \in \mathbb{V}$ as

$$(3.3) \qquad \pi[z, y] \doteq \begin{cases} 0, & \text{if } z_i = 0 \text{ and } y = v_i \text{ for some } i = 1, 2, \\ y, & \text{otherwise.} \end{cases}$$

The distribution of $Z$ is completely determined by that of the sequence $Y = \{Y(k)\}$. Define

$$\mathcal{P}^+(\mathbb{V}) \doteq \{\theta = (\theta_0, \theta_1, \theta_2) : \theta \text{ is a probability measure on } \mathbb{V}$$
$$\text{and } \theta_i = \theta[v_i] > 0 \text{ for every } i = 0, 1, 2\}.$$

Under the (true) probability measure $\mathbb{P}$, $Y$ is a sequence of independent identically distributed (i.i.d.) random variables with distribution

$$\Theta \doteq (\lambda, \mu_1, \mu_2) \in \mathcal{P}^+(\mathbb{V}).$$

3.3. *The dynamic importance sampling algorithms.* The importance sampling schemes we consider use state-dependent changes of measure that can be characterized by stochastic kernels $\bar{\Theta}^n[\cdot|\cdot]$ on $\mathbb{V}$ given $\mathbb{R}^2_+$, that is, $\bar{\Theta}^n[\cdot|x] \in \mathcal{P}^+(\mathbb{V})$ for every $x \in \mathbb{R}^2_+$.

To be more precise, for a given threshold $n$, define the scaled state process $X^n = Z/n$, where $Z$ is defined as in (3.2). Since the definition of $\pi$ implies $\pi[nx, y] = \pi[x, y]$ for every $x \in \mathbb{R}^2_+$, $X^n$ satisfies the equation

$$(3.4) \qquad X^n(k+1) = X^n(k) + \frac{1}{n}\pi[X^n(k), Y(k+1)],$$

with initial condition $X^n(0) = Z(0)/n = 0$. The importance sampling generates $\{Y(k)\}$ as follows. The conditional probability of $Y(k+1) = v_i$, given $\{Y(j) : j = 1, 2, \ldots, k\}$, is just $\bar{\Theta}^n[v_i|X^n(k)]$ for each $i = 0, 1, 2$.

Define the hitting times

$$T_n \doteq \inf\{k \geq 0 : X_1^n(k) + X_2^n(k) = 1\},$$
$$T_0 \doteq \inf\{k \geq 1 : X_1^n(k) = X_2^n(k) = 0\}.$$

Let $A_n$ be the event of interest, that is,

$$A_n = \{X_1^n + X_2^n \text{ reaches } 1 \text{ before returning to } 0\} = \{T_n < T_0\}.$$

The importance sampling estimator is just

$$(3.5) \qquad \hat{p}_n = 1_{A_n} \cdot \prod_{k=0}^{T_n - 1} \frac{\Theta[Y(k+1)]}{\bar{\Theta}^n[Y(k+1)|X^n(k)]}.$$

The second moment of $\hat{p}_n$, thanks to (2.1), equals $E^{\mathbb{P}}[\hat{p}_n]$. The goal is to choose a stochastic kernel $\bar{\Theta}^n$ so that this second moment (whence the variance of $\hat{p}_n$) is as small as possible. Another important consideration is that one would like $\bar{\Theta}^n$ to be simple and easy to implement.

3.4. *Notation and terminology.* Before we proceed to construct importance sampling algorithms, we collect in this section some notation and terminology. Define

$$\bar{D} \doteq \{(x_1, x_2) : x_i \geq 0, x_1 + x_2 \leq 1\},$$
$$D \doteq \{(x_1, x_2) : x_i > 0, x_1 + x_2 < 1\},$$
$$\partial_1 \doteq \{(0, x_2) : 0 < x_2 < 1\},$$
$$\partial_2 \doteq \{(x_1, 0) : 0 < x_1 < 1\},$$
$$\partial_e \doteq \{(x_1, x_2) : x_i \geq 0, x_1 + x_2 = 1\},$$
$$\bar{D}_n \doteq \bar{D} \cap \{(z_1, z_2)/n : (z_1, z_2) \in \mathbb{Z}^2_+\},$$
$$D_n \doteq D \cap \{(z_1, z_2)/n : (z_1, z_2) \in \mathbb{Z}^2_+\}.$$

Sometimes we refer to $\partial_e$ as the "exit boundary."

REMARK 3.1. *Relative entropy representation for exponential integrals.* Let $(S, \mathcal{F})$ be a measurable space and $f: S \to \mathbb{R}$ a bounded measurable function. Denote by $\mathcal{P}(S)$ the space of probability measures on $(S, \mathcal{F})$. Then for any $\gamma \in \mathcal{P}(S)$,

$$-\log \int_S e^{-f} \, d\gamma = \inf_{\theta \in \mathcal{P}(S)} \left[ R(\theta \| \gamma) + \int_S f \, d\theta \right].$$

Furthermore, the minimizer of the right-hand side exists and is mutually absolutely continuous with respect to $\gamma$. Here the relative entropy $R(\cdot \| \cdot)$ is defined as

$$R(\theta \| \gamma) \doteq \begin{cases} \int_S \log \frac{d\theta}{d\gamma} \, d\theta, & \text{if } \theta \ll \gamma, \\ \infty, & \text{otherwise.} \end{cases}$$

We refer the readers to [4], Proposition 1.4.2, for the proof.

3.5. *The Isaacs equation.* In this section we formally derive the Isaacs equation associated with the limit differential game that lies underneath importance sampling algorithms. A rigorous argument, though possible, is not necessary for our purpose.

Recall our goal is to choose a stochastic kernel $\bar{\Theta}^n$ so as to keep the second moment $E^{\mathbb{P}}[\hat{p}_n]$ as small as possible. We can think of this as a stochastic control problem and write down the corresponding dynamic programming equation (DPE). To this end, we extend the dynamics and let, for every $x \in \bar{D}_n$,

$$V_n(x) \doteq \inf_{\bar{\Theta}^n} E_x^{\mathbb{P}}[\hat{p}_n] = \inf_{\bar{\Theta}^n} E_x^{\mathbb{P}}\left[ 1_{A_n} \cdot \prod_{k=0}^{T_n - 1} \frac{\Theta[Y(k+1)]}{\bar{\Theta}^n[Y(k+1) | X^n(k)]} \right],$$

where $\hat{p}_n$ is defined in exactly the same fashion as in Section 3.3 and $E_x^{\mathbb{P}}$ denotes expected value conditioned on $X^n(0) = x$.

For simplicity, we further assume that $x \in D_n$, whence $\pi[x, y] \equiv y$ for every $y \in \mathbb{V}$. Under the original probability measure $\mathbb{P}$, the sequence $\{Y(k)\}$ is i.i.d. with distribution $\Theta$. Hence the DPE

$$V_n(x) = \inf_{\bar{\Theta} \in \mathcal{P}^+(\mathbb{V})} \sum_{i=0}^{2} V_n\left(x + \frac{1}{n} v_i\right) \frac{\Theta[v_i]}{\bar{\Theta}[v_i]} \cdot \Theta[v_i]$$

holds. Consider a logarithmic transform of $V_n$ and define

$$W_n(x) \doteq -\frac{1}{n} \log V_n(x).$$

We have

$$W_n(x) = \sup_{\bar{\Theta} \in \mathcal{P}^+(\mathbb{V})} -\frac{1}{n} \log \sum_{i=0}^{2} \exp\left\{ -n W_n\left(x + \frac{1}{n} v_i\right) - \log \frac{\bar{\Theta}[v_i]}{\Theta[v_i]} \right\} \Theta[v_i].$$

Applying the relative entropy representation for exponential integrals (see Remark 3.1) to the right-hand side of the last equation, it follows that

$$W_n(x) = \sup_{\bar{\Theta} \in \mathcal{P}^+(\mathbb{V})} \inf_{\theta \in \mathcal{P}^+(\mathbb{V})} \left[ \sum_{i=0}^2 W_n\left(x + \frac{1}{n}v_i\right)\theta[v_i] \right.$$
$$\left. + \frac{1}{n}\left(\sum_{i=0}^2 \theta[v_i]\log\frac{\bar{\Theta}[v_i]}{\Theta[v_i]} + R(\theta\|\Theta)\right)\right].$$

Note that taking infimum over $\theta \in \mathcal{P}^+(\mathbb{V})$ is equivalent to taking infimum over $\theta \in \mathcal{P}(\mathbb{V})$ since by Remark 3.1 the minimizing $\theta$ is mutually absolutely continuous to $\Theta$, whence it belongs to $\mathcal{P}^+(\mathbb{V})$.

Suppose for now that $W_n(x)$ converges to $W(x)$. Let $DW$ be the gradient of $W$ and formally assume the approximation

$$W_n\left(x + \frac{1}{n}v_i\right) - W_n(x) \approx \frac{1}{n}\langle DW(x), v_i\rangle.$$

Observing $\sum \theta[v_i] = 1$, we arrive at

$$(3.6) \quad 0 = \sup_{\bar{\Theta} \in \mathcal{P}^+(\mathbb{V})} \inf_{\theta \in \mathcal{P}^+(\mathbb{V})} \left[\langle DW(x), \mathbb{F}(\theta)\rangle + \sum_{i=0}^2 \theta[v_i]\log\frac{\bar{\Theta}[v_i]}{\Theta[v_i]} + R(\theta\|\Theta)\right],$$

where

$$(3.7) \quad \mathbb{F}(\theta) \doteq \sum_{i=0}^2 \theta[v_i] \cdot v_i$$

for each $\theta \in \mathcal{P}^+(\mathbb{V})$. Equation (3.6) is called an *Isaacs equation*.

We now discuss the boundary conditions. For the exit boundary, we have by definition $V_n(x) = 1$ or $W_n(x) = 0$, therefore we impose the Dirichlet boundary condition

$$(3.8) \quad W(x) = 0 \quad \text{for } x \in \partial_\mathsf{e}.$$

For $\partial_1$ and $\partial_2$, we impose the Neumann boundary condition that is typically associated with constrained dynamics [12]

$$(3.9) \quad \langle DW(x), d_i\rangle = 0 \quad \text{for } x \in \partial_i.$$

Finally, we make a few remarks on the game interpretation of importance sampling. The Isaacs equation (3.6) indicates that the underlying game has two players. The player who chooses the change of measure in order to minimize the second moment (i.e., $\bar{\Theta}$) becomes the maximizing player in the game due to the negative sign in the logarithmic transform. The minimizing player is artificially introduced, and chooses $\theta$. We will refer to this player as the "large deviation player." The dynamics of the game are completely determined by $\theta$, or the choice of the large deviation player, while the running cost of the game depends on the choices of both players.

3.6. *The properties of the Hamiltonian.* Our construction of importance sampling algorithms is based on classical subsolutions to the Isaacs equation. Therefore it is useful to study the properties of this equation. Define for each $p \in \mathbb{R}^2$

$$(3.10) \quad \mathbb{H}(p) \doteq \sup_{\bar{\Theta} \in \mathcal{P}^+(\mathbb{V})} \inf_{\theta \in \mathcal{P}^+(\mathbb{V})} \left[ \langle p, \mathbb{F}(\theta) \rangle + \sum_{i=0}^{2} \theta[v_i] \log \frac{\bar{\Theta}[v_i]}{\Theta[v_i]} + R(\theta \| \Theta) \right].$$

The function $\mathbb{H}$ is called the Hamiltonian, and the Isaacs equation (3.6) can be written as

$$(3.11) \quad \mathbb{H}(DW) = 0.$$

We have the following result, whose proof is straightforward and therefore omitted.

PROPOSITION 3.2. *Let $\mathbb{H}$ be defined as in (3.10).*

1. *For each $p = (p_1, p_2) \in \mathbb{R}^2$, there exists a saddle point $(\bar{\Theta}^*(p), \theta^*(p)) \in \mathcal{P}^+(\mathbb{V}) \times \mathcal{P}^+(\mathbb{V})$ given by*

$$\bar{\Theta}^*(p) = \theta^*(p) = N(p)(\lambda e^{-p_1/2}, \mu_1 e^{(p_1-p_2)/2}, \mu_2 e^{p_2/2}),$$

*where*

$$N(p) \doteq [\lambda e^{-p_1/2} + \mu_1 e^{(p_1-p_2)/2} + \mu_2 e^{p_2/2}]^{-1}.$$

*In particular, the order of sup and inf can be exchanged in (3.10).*
2. $\mathbb{H}$ *is concave and has the representation*

$$\mathbb{H}(p) = \inf_{\theta \in \mathcal{P}^+(\mathbb{V})} [\langle p, \mathbb{F}(\theta) \rangle + 2R(\theta \| \Theta)] = 2 \log N(p).$$

For any $p \in \mathbb{R}^2$, we will refer to $\bar{\Theta}^*(p)$ as the (saddle point) change of measure corresponding to $p$.

Figure 3 is a picture of the zero-level curve of $\mathbb{H}$. Recall that $\gamma$, as defined in (3.1), equals $\log(\mu_2/\lambda)$.

3.7. *The solution to the Isaacs equation.* It is well known that viscosity solutions provide physically meaningful solutions to equations such as (3.11). However, viscosity solutions to the Isaacs equation (3.11) and boundary conditions (3.8) and (3.9), which are only weak-sense solutions, are *not* suitable for the purpose of constructing efficient importance sampling algorithms for this tandem Jackson network.

More precisely, consider the very simple, affine function

$$W_{\mathsf{s}}(x) \doteq \langle r_1, x \rangle + 2\gamma.$$

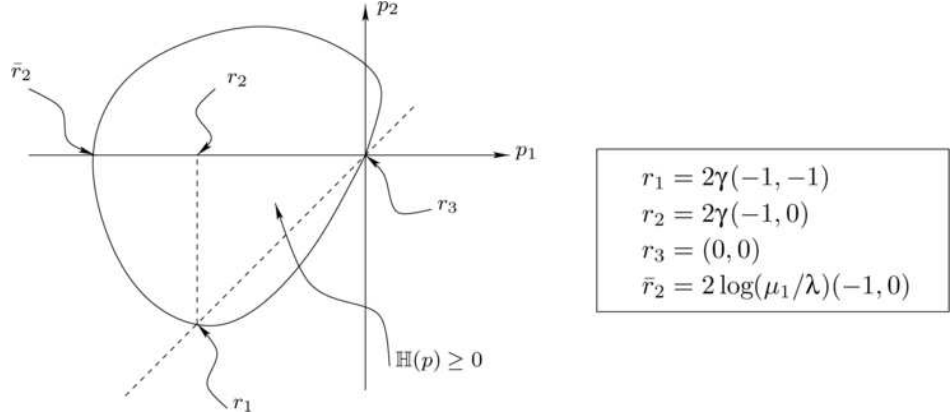

Fig. 3. *Hamiltonian* $\mathbb{H}$.

$W_{\mathsf{s}}$ is a viscosity solution to the Isaacs equation (3.11) and boundary conditions (3.8) and (3.9). Even though $W_{\mathsf{s}}(0)$ equals the optimal decay rate $2\gamma$, the corresponding saddle point change of measure, by Proposition 3.2, is

$$\bar{\Theta}^*(DW_{\mathsf{s}}) = \bar{\Theta}^*(r_1) = (\mu_2, \mu_1, \lambda),$$

which is exactly the state-independent change of measure based on standard heuristic (i.e., switching the arrival rate and the smallest service rate) and therefore inefficient in general.

As remarked previously, the failure of the importance sampling based on $W_{\mathsf{s}}$ is due to the fact that $W_{\mathsf{s}}$ is only a weak-sense viscosity solution. It is not a classical solution (or even a classical subsolution as defined in the next subsection), since on the boundary $\partial_2$

$$\langle DW_{\mathsf{s}}, d_2 \rangle = \langle r_1, d_2 \rangle = -2\gamma < 0.$$

In a sense that we will make precise later on, this inequality is in the "wrong" direction, which suggests that the (artificial) large deviation player, who determines the game dynamics, may be able to exploit this "bad" boundary to a degree that the importance sampling estimator based on $W_{\mathsf{s}}$ becomes inefficient. It is not coincidental, as observed in [10], that the inefficiency is because a sample path can spend a significant amount of time near the boundary $\partial_2$ before leaving domain $D$ and thereby accumulate a huge Radon–Nikodým derivative.

3.8. *Subsolutions and importance sampling schemes.* The idea of [9] is that classical subsolutions to Isaacs equations can be used to construct efficient importance sampling schemes. It has advantages over solution-based importance sampling schemes in simplicity, greater flexibility, and general applicability. The goal of this section is to construct classical subsolutions and identify the corresponding changes of measure. As in [9], the construction is divided into two steps. We first identify a piecewise smooth subsolution as the minimum of affine functions and then mollify it to obtain a classical subsolution.

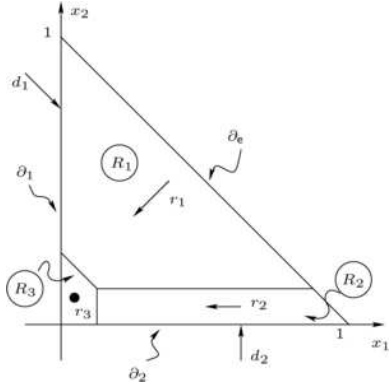

Fig. 4. *Piecewise affine subsolution.*

DEFINITION 3.3. A function $W : \bar{D} \to \mathbb{R}$ is a classical subsolution to the Isaacs equation (3.11) and boundary conditions (3.8) and (3.9) if:

1. $W$ is continuously differentiable,
2. $\mathbb{H}(DW(x)) \geq 0$ for every $x \in D$,
3. $W(x) \leq 0$ for $x \in \partial_{\mathsf{e}}$,
4. $\langle DW(x), d_i \rangle \geq 0$ for $x \in \partial_i$, $i = 1, 2$.

3.8.1. *Construction of piecewise affine subsolutions.* We will need a piecewise affine subsolution $\bar{W}$ with the following properties (see Figure 4).

1. $\bar{W}$ can be written as $\bar{W} = \bar{W}_1 \wedge \bar{W}_2 \wedge \bar{W}_3$ where $\bar{W}_k$ is an affine function for each $k = 1, 2, 3$.
2. $\bar{D}$ is divided into three regions $R_1$, $R_2$ and $R_3$, such that in each region $R_k$, $\bar{W} = \bar{W}_k$.
3. The subsolution property $\mathbb{H}(D\bar{W}(x)) = \mathbb{H}(D\bar{W}_k(x)) \geq 0$ holds for every $x$ in the interior of region $R_k$.
4. The Dirichlet boundary inequality $\bar{W}(x) \leq 0$ holds for $x \in \partial_{\mathsf{e}}$.
5. The Neumann boundary inequality $\langle D\bar{W}(x), d_i \rangle \geq 0$ holds whenever $x \in \partial_i$ and $D\bar{W}(x)$ exists.

One such subsolution can be constructed as follows. Fix an arbitrary $\delta > 0$ and let, for each $k$,

$$\bar{W}_k^\delta(x) \doteq \langle r_k, x \rangle + 2\gamma - k\delta, \tag{3.12}$$

where the $r_k$'s are depicted in Figure 3. It is not difficult to check that
$$\bar{W}^\delta \doteq \bar{W}_1^\delta \wedge \bar{W}_2^\delta \wedge \bar{W}_3^\delta$$
satisfies all the requirements for all small $\delta > 0$.

REMARK 3.4. The failure of the boundary inequality along the $x_1$ axis for $W_{\mathsf{s}}$ (see Section 3.7) corresponds to the existence of a boundary layer in the prelimit which vanishes in the limit. It is for this reason that we introduce $\bar{W}_2^\delta$, which perturbs the gradient in a neighborhood of this axis. A similar perturbation is not required along the $x_2$ axis, since the boundary inequality already holds there. $\bar{W}_3^\delta$ is introduced to ensure that both boundary conditions hold in a neighborhood of the origin.

3.8.2. *Mollification.* To mollify the piecewise affine subsolution $\bar{W}^\delta$, we will adopt a mollification called *exponential weighting* that is specialized to the minimum of a finite set of smooth functions. For future reference, we describe the mollification in its general form.

Consider continuously differentiable functions $\{h_1, h_2, \ldots, h_K\}$, and let
$$h \doteq h_1 \wedge h_2 \wedge \cdots \wedge h_K.$$
Fix a small positive number $\varepsilon$ and define the mollification
$$h^\varepsilon(x) \doteq -\varepsilon \log \sum_{k=1}^{K} \exp\left\{-\frac{1}{\varepsilon} h_k(x)\right\}.$$
We have the following result, whose proof is straightforward and can be found in [9], Section 3.3.

LEMMA 3.5. *For any $\varepsilon > 0$, $h^\varepsilon$ is continuously differentiable with*
$$Dh^\varepsilon(x) = \sum_{k=1}^{K} \rho_k^\varepsilon(x) D h_k(x),$$
*where*
$$\rho_i^\varepsilon(x) \doteq \frac{\exp\{-h_i(x)/\varepsilon\}}{\sum_{k=1}^{K} \exp\{-h_k(x)/\varepsilon\}}.$$
*Furthermore, we have the uniform bounds*
$$-K\varepsilon \leq h^\varepsilon(x) - h(x) \leq 0$$
*for every $x$.*

Note that $\rho^\varepsilon(x) \doteq (\rho_1^\varepsilon(x), \rho_2^\varepsilon(x), \ldots, \rho_K^\varepsilon(x))$ defines a probability vector in the sense that $\rho_k^\varepsilon(x) \geq 0$ and
$$\sum_{k=1}^{K} \rho_k^\varepsilon(x) \equiv 1.$$

3.8.3. *The classical subsolution.* Applying this mollification to $\bar{W}^\delta$, we define

$$(3.13) \qquad W^{\varepsilon,\delta}(x) \doteq -\varepsilon \log \sum_{k=1}^{3} \exp\left\{-\frac{1}{\varepsilon}\bar{W}_k^\delta(x)\right\}.$$

Thanks to Lemma 3.5, $W^{\varepsilon,\delta}$ is continuously differentiable with

$$(3.14) \qquad DW^{\varepsilon,\delta}(x) = \sum_{k=1}^{3} \rho_k^{\varepsilon,\delta}(x) r_k,$$

where

$$(3.15) \qquad \rho_i^{\varepsilon,\delta}(x) \doteq \frac{\exp\{-\bar{W}_i^\delta(x)/\varepsilon\}}{\sum_{k=1}^{3} \exp\{-\bar{W}_k^\delta(x)/\varepsilon\}}.$$

We should notice that with this mollification, the function $\bar{W}^{\varepsilon,\delta}$ is not precisely a classical subsolution, but only approximately. Indeed, Lemma B.1 states that the Neumann boundary conditions are only satisfied approximately in the sense that, for $x \in \partial_i$,

$$\langle DW^{\varepsilon,\delta}(x), d_i \rangle \geq -\bar{\varepsilon}$$

for some small positive number $\bar{\varepsilon}$ as long as $\varepsilon/\delta$ is chosen small. The reason for this violation of the subsolution property is that the exponential weighting is not a "local" smoothing. However, the advantages of the exponential weighting (especially the analytical tractability) outweigh the minor additional complications in the analysis introduced by this error.

3.8.4. *The importance sampling estimator and its asymptotics.* For each $k$, let $\bar{\Theta}_k^*$ be the saddle point change of measure that corresponds to the affine function $\bar{W}_k$, that is,

$$\bar{\Theta}_k^* \doteq \bar{\Theta}^*(D\bar{W}_k) = \bar{\Theta}^*(r_k) \in \mathcal{P}^+(\mathbb{V}),$$

where $\bar{\Theta}^*(\cdot)$ is as defined in Proposition 3.2. Straightforward calculation yields that

$$\bar{\Theta}_1^* = (\mu_2, \mu_1, \lambda), \qquad \bar{\Theta}_2^* = \frac{1}{\lambda\mu_1 + 2\mu_2^2}(\mu_2^2, \lambda\mu_1, \mu_2^2), \qquad \bar{\Theta}_3^* = (\lambda, \mu_1, \mu_2).$$

The change of measure based on the $W^{\varepsilon,\delta}$ is just a state-dependent mixture of $\bar{\Theta}_k^*$. More precisely, define a stochastic kernel $\bar{\Theta}^{\varepsilon,\delta}[\cdot|\cdot]$ by

$$(3.16) \qquad \bar{\Theta}^{\varepsilon,\delta}[\cdot|x] \doteq \sum_{k=1}^{3} \rho_k^{\varepsilon,\delta}(x) \bar{\Theta}_k^*(\cdot) \in \mathcal{P}^+(\mathbb{V}),$$

and for each fixed $n$, let

$$\bar{\Theta}^n[\cdot|\cdot] \equiv \bar{\Theta}^{\varepsilon,\delta}[\cdot|\cdot]. \tag{3.17}$$

In other words, the importance sampling algorithm simulates $Y(k+1)$, conditional on the sample history $\{Y(j): 1 \leq j \leq k\}$, from the distribution $\bar{\Theta}^{\varepsilon,\delta}[\cdot|X^n(k)]$, where $X^n$ is the state process as defined in (3.4). The importance sampling estimator $\hat{p}_n$ is then given by (3.5).

We have the following result regarding its asymptotic performance, whose proof is deferred to Appendix B.

THEOREM 3.6. *There exist a pair of positive constants $(A, B)$ that only depend on the system parameters $(\lambda, \mu_1, \mu_2)$ such that, provided $\varepsilon/\delta < B$, the second moment of the importance sampling estimator $\hat{p}_n$ satisfies*

$$\liminf_n -\frac{1}{n}\log[\textit{2nd moment of } \hat{p}_n] \geq 2\gamma - F(\varepsilon, \delta),$$

*where*

$$F(\varepsilon, \delta) \doteq 3\varepsilon + 3\delta + A\exp\{-\delta/\varepsilon\}.$$

Since $2\gamma$ is the optimal decay rate, the theorem suggest that the importance sampling scheme is nearly asymptotically optimal as long as $F(\varepsilon, \delta)$ is small. This can be achieved if one sets both $\delta$ and $\varepsilon/\delta$ small.

REMARK 3.7. The formula of $F$ also provides an interesting relation between $\varepsilon$ and $\delta$. For each fixed small $\varepsilon$, $F(\varepsilon, \cdot)$ is minimized at

$$\delta = -\varepsilon \log \varepsilon + \varepsilon \log \frac{A}{3} \approx -\varepsilon \log \varepsilon.$$

This suggests that a good strategy is to set $\delta = -\varepsilon \log \varepsilon$. Note that in this case, when $\varepsilon$ is small, so are $\delta$ and $F(\varepsilon, \delta)$.

3.8.5. *Asymptotic optimality.* The previous section provides a nearly asymptotically optimal importance sampling algorithm. It is good enough for many practical purposes where $n$ is large but not exceedingly large. However, one would still like to see an algorithm that gives optimality. This only requires a slight modification.

Instead of using a fixed pair of parameters $\varepsilon$ and $\delta$ for all $n$, we now allow them to vary depending on $n$ and denote them by $\varepsilon_n$ and $\delta_n$. For each $n$, we use the change of measure based on $W^{\varepsilon_n, \delta_n}$, which amounts to letting

$$\bar{\Theta}^n[\cdot|x] \doteq \bar{\Theta}^{\varepsilon_n,\delta_n}[\cdot|x] = \sum_{k=1}^{3} \rho_k^{\varepsilon_n,\delta_n}(x) \bar{\Theta}_k^*(\cdot) \in \mathcal{P}^+(\mathbb{V}). \tag{3.18}$$

Abusing the notation a bit, we again denote by $\hat{p}_n$ the corresponding importance sampling estimator.

THEOREM 3.8. *The estimator $\hat{p}_n$ is asymptotically optimal, that is,*

$$\lim_n -\frac{1}{n}\log[\text{2nd moment of } \hat{p}_n] = 2\gamma,$$

*provided that $\delta_n \to 0$, $\varepsilon_n/\delta_n \to 0$ and $n\varepsilon_n \to \infty$.*

Remark 3.7 suggests that a good choice is to set $\delta_n = -\varepsilon_n \log \varepsilon_n$. In this case, asymptotic optimality follows if $\varepsilon_n \to 0$ and $n\varepsilon_n \to \infty$.

3.8.6. *Further remarks on the importance sampling algorithms.* The computation of the weights $\{\rho_k^{\varepsilon,\delta}\}$ or $\{\rho_k^{\varepsilon_n,\delta_n}\}$ is very simple. As a consequence, the dynamic importance sampling algorithms based on (3.16), (3.17) or (3.18) are practically as fast as the standard heuristic scheme where a constant change of measure is used.

It is possible that one can associate other changes of measure with subsolutions. For example, one can define $\bar{\Theta}^n[\cdot|x] \equiv \bar{\Theta}^*(DW^{\varepsilon,\delta}(x))$ in lieu of (3.16) and (3.17), and the resulting algorithms will have similar asymptotic performance. However, the use of mixtures such as (3.16) is computationally more convenient. This is especially the case when the change of measure that corresponds to a particular gradient is not easily obtainable. For example, for a system with Markov modulated arrival and service rates, the computation of the change of measure corresponding to a single gradient $p$ requires solving an eigenvalue/eigenvector problem. If we smooth first and then compute the change of measure suitable for each point $x$, then many such problems must be solved. In contrast, mixtures like (3.16) only require the computation of the changes of measure that correspond to the finite collection of vectors $r_k$.

3.9. *Numerical results.* In this section we present some numerical results for the case where $\lambda = 0.1$, $\mu_1 = \mu_2 = 0.45$. For comparison, the theoretical value of $p_n$ is obtained by iteratively solving the linear system of equations that characterize this probability, an approach that is feasible when the system is sufficiently small. Note that in this case, the standard heuristic importance sampling scheme leads to estimators with infinite variance [3].

In the simulations, we always set $\delta = -\varepsilon \log \varepsilon$. This choice was suggested by Remark 3.7, and was experimentally observed to be a good choice for small $\varepsilon$. We ran simulations for $n = 20$, with $\varepsilon = 0.01, 0.02$ and $0.03$, respectively. For each $\varepsilon$ we present two estimates and each estimate consists of 20,000 replications. The theoretical is $p_n = 6.0 \times 10^{-12}$ (see Table 1).

In all the tables, "Std. Err." stands for "standard error" and "C.I." for "confidence interval." The performance of the dynamic importance sampling schemes based on subsolutions is stable across different simulations, with good estimates and small standard errors.

In Table 2 there are more simulation results with $n = 30, 40, 50$, with $\varepsilon = 0.02$ and $\delta = -\varepsilon \log \varepsilon$. Each estimate consists of 20,000 replications.

REMARK 3.9. It is not difficult to check that the "thickness" or the height of the boundary region $R_2$ in Figure 4 is $\delta/(2\gamma)$. Since we are scaling the queue sizes by a factor $n$, the thickness of the boundary region in the prelimit will be $n\delta/(2\gamma)$ when unscaled. However, the optimality conditions $n\varepsilon_n \to \infty$ and $\varepsilon_n/\delta_n \to 0$ in Theorem 3.8 imply that $n\delta_n \to \infty$. This does not allow the boundary region to be too thin in the prelimit. The need for such control is supported by experimentation, which shows that for a fixed $n$, the simulation results tend to deteriorate when $\varepsilon$ is too small.

**4. Extensions to $d$-node tandem Jackson networks.** The work on the two-node tandem Jackson network can be easily extended to $d$-node tandem Jackson networks and more general exit probabilities. To be more precise, consider a $d$-node tandem Jackson network with Poisson arrival rate $\lambda$ and consecutive exponential service rates $\mu_1, \ldots, \mu_d$. The state of the network is described by the embedded Markov chain $Z = \{Z(k)\} = \{(Z_1(k), \ldots, Z_d(k))\}$, where $Z_i$ denotes the queue length at node $i$. The system is assumed to be stable, that is, $\lambda < \min\{\mu_1, \ldots, \mu_d\}$. Let $\Gamma$ be a closed subset of $\mathbb{R}_+^d$ such that $0 \notin \Gamma$ and the closure of $\mathbb{R}_+^d \setminus \Gamma$ is compact. We are interested in the following rare-event probability:

$$p_n \doteq \mathbb{P}\{\text{Process } Z \text{ hits set } n\Gamma \text{ before returning to } 0, \text{ starting from } 0\}.$$

TABLE 1
*IS based on subsolutions, two-node tandem, total population overflow*

|  | $\varepsilon = 0.01$ | | $\varepsilon = 0.02$ | | $\varepsilon = 0.03$ | |
| --- | --- | --- | --- | --- | --- | --- |
|  | No. 1 | No. 2 | No. 1 | No. 2 | No. 1 | No. 2 |
| Estimate ($\times 10^{-12}$) | 5.7 | 5.5 | 5.8 | 6.1 | 6.3 | 5.8 |
| Std. Err. ($\times 10^{-12}$) | 0.4 | 0.3 | 0.3 | 0.5 | 0.4 | 0.2 |
| 95% C.I. ($\times 10^{-12}$) | [4.9, 6.4] | [4.9, 6.1] | [5.2, 6.4] | [5.1, 7.1] | [5.5, 7.1] | [5.3, 6.3] |

TABLE 2
*IS based on subsolutions, two-node tandem, total population overflow*

|  | $n = 30$ | $n = 40$ | $n = 50$ |
| --- | --- | --- | --- |
| Theoretical value | $2.63 \times 10^{-18}$ | $1.03 \times 10^{-24}$ | $3.80 \times 10^{-31}$ |
| Estimate | $2.73 \times 10^{-18}$ | $1.05 \times 10^{-24}$ | $3.75 \times 10^{-31}$ |
| Std. Err. | $0.18 \times 10^{-18}$ | $0.03 \times 10^{-24}$ | $0.16 \times 10^{-31}$ |
| 95% C.I. | $[2.37, 3.09] \times 10^{-18}$ | $[0.99, 1.11] \times 10^{-24}$ | $[3.43, 4.07] \times 10^{-31}$ |

Without loss of generality, we assume that $\lambda + \mu_1 + \cdots + \mu_d = 1$. We also assume that $p_n$ decays exponentially with

$$\lim_n -\frac{1}{n} \log p_n = \gamma.$$

4.1. *Isaacs equation and the Hamiltonian.* The increments of $Z$ take values in $\mathbb{V} = \{v_0, v_1, \ldots, v_d\}$ where the $v_i$'s are $d$-dimensional vectors defined by

$$v_0 = (1, 0, \ldots, 0), \qquad [v_i]_j \doteq \begin{cases} -1, & \text{if } j = i, \\ 1, & \text{if } j = i+1 \text{ and } j \leq d, \\ 0, & \text{otherwise.} \end{cases}$$

Similarly to (3.4), the scaled state process $X^n \doteq Z/n$ satisfies

$$X^n(k+1) = X^n(k) + \frac{1}{n} \pi[X^n(k), Y(k+1)],$$

where $\pi$ plays the same role as in (3.3). The sequence $\{Y(k)\}$ consists of i.i.d. random variables taking values in $\mathbb{V}$ with common distribution

$$\Theta = (\lambda, \mu_1, \ldots, \mu_d) \in \mathcal{P}^+(\mathbb{V}).$$

Define the regions

$$D \doteq \{x \in \mathbb{R}_+^d : x \notin \Gamma, x_i > 0, i = 1, \ldots, d\},$$
$$\partial_i \doteq \{x \in \mathbb{R}_+^d : x \notin \Gamma, x_i = 0\}, i = 1, \ldots, d,$$

and the directions of constraints $d_i = -v_i$.

The Isaacs equation is just $\mathbb{H}(DW) = 0$, where

$$\mathbb{H}(p) = \sup_{\bar{\Theta} \in \mathcal{P}^+(\mathbb{V})} \inf_{\theta \in \mathcal{P}^+(\mathbb{V})} \left[ \langle p, \mathbb{F}(\theta) \rangle + \sum_{i=0}^d \theta[v_i] \log \frac{\bar{\Theta}[v_i]}{\Theta[v_i]} + R(\theta \| \Theta) \right],$$

with

$$\mathbb{F}(\theta) \doteq \sum_{i=0}^d \theta[v_i] \cdot v_i$$

for every $\theta \in \mathcal{P}^+(\mathbb{V})$. The boundary conditions are $W(x) = 0$ for $x \in \Gamma$ and $\langle DW(x), d_i \rangle = 0$ for $x \in \partial_i$.

The following result is an extension of Proposition 3.2, whose proof is very similar and thus omitted.

PROPOSITION 4.1. *For every $p \in \mathbb{R}^d$, there exists a saddle point for the Hamiltonian $\mathbb{H}$, say $(\bar{\Theta}^*(p), \theta^*(p)) \in \mathcal{P}^+(\mathbb{V}) \times \mathcal{P}^+(\mathbb{V})$, given by*

$$\bar{\Theta}^*(p)[v_i] = \theta^*(p)[v_i] = N(p) \cdot \Theta[v_i] \exp\{-\langle p, v_i \rangle / 2\},$$

*where*

$$N(p) \doteq \left[\sum_{i=0}^{d} \Theta[v_i] \exp\{-\langle p, v_i\rangle/2\}\right]^{-1}.$$

*Moreover, the Hamiltonian $\mathbb{H}$ is concave and $\mathbb{H}(p) = 2\log N(p)$.*

4.2. *Subsolutions and importance sampling schemes.* The construction of subsolution proceeds in an analogous fashion: we start with a piecewise smooth subsolution and then mollify it by exponential weighting. We will discuss the general case where the subsolutions can vary depending on $n$. To be more specific, let $(\bar{W}_1^n, \ldots, \bar{W}_K^n)$ be smooth functions (preferably affine functions) and let

$$\bar{W}^n \doteq \bar{W}_1^n \wedge \cdots \wedge \bar{W}_K^n.$$

The choice of $\{\bar{W}_k^n\}$ should have the following properties:

1. $\mathbb{H}(D\bar{W}_k^n(x)) \geq 0$ for every $x \in D$ and every $k$,
2. $\bar{W}^n(x) \leq 0$ for every $x \in \Gamma$,
3. for $x$ on boundary $\partial_i$, $\langle D\bar{W}^n(x), d_i\rangle \geq 0$ when $D\bar{W}^n(x)$ is well defined.

The quantities $W^{\varepsilon_n,n}(x), \rho_i^{\varepsilon_n,n}(x)$ and the stochastic kernel $\bar{\Theta}^n$ are defined in a fashion exactly analogous to (3.13), (3.15), (3.16) and (3.17). We denote by $\hat{p}_n$ the corresponding importance sampling estimator.

The following result is an extension of Theorem 3.8. The proof is very similar and thus omitted.

THEOREM 4.2. *Assume that $\{\bar{W}_k^n(x)\}$ has uniformly bounded first and second derivatives for $x \in D$ and that there exists $\bar{\varepsilon}_n \geq 0$ such that for $x \in \partial_i$, $\langle DW^{\varepsilon_n,n}(x), d_i\rangle \geq -\bar{\varepsilon}_n$. We also assume that $\liminf_n \bar{W}^n(0) \geq 2\gamma$. Then the importance sampling estimator $\hat{p}_n$ is asymptotically optimal, that is,*

$$\lim_n -\frac{1}{n}\log[\text{2nd moment of } \hat{p}_n] = 2\gamma,$$

*provided that $\varepsilon_n \to 0$, $\bar{\varepsilon}_n \to 0$ and $n\varepsilon_n \to \infty$.*

REMARK 4.3. One can also write down a result similar to Theorem 3.6 for the case where $\bar{W}_k^n \equiv \bar{W}_k$ and $\varepsilon_n \equiv \varepsilon$, $\bar{\varepsilon}_n \equiv \bar{\varepsilon}$. The corresponding importance sampling estimator, still denoted by $\hat{p}_n$, will satisfy

$$\liminf_n -\frac{1}{n}\log[\text{2nd moment of } \hat{p}_n] \geq \bar{W}(0) - (K\varepsilon + C\bar{\varepsilon}),$$

where $C$ is a constant only depends on the system parameter $\Theta$, under the condition that $\bar{\varepsilon}$ is small enough.

TABLE 3

|  | $r_k$ | $\bar{\Theta}^*(r_k)$ |
|---|---|---|
| $\mu_1 \geq \mu_2$ | $r_1 = 2\log(\mu_2/\lambda)(-1,-1)$<br>$r_2 = 2\log(\mu_1/\lambda)(-1,0)$<br>$r_3 = (0,0)$ | $(\mu_2, \mu_1, \lambda)$<br>$(\mu_1, \lambda, \mu_2)$<br>$(\lambda, \mu_1, \mu_2)$ |
| $\mu_1 < \mu_2$ | $r_1 = (-2\log(\mu_1/\lambda), -2\log(\mu_2/\lambda))$<br>$r_2 = 2\log(\mu_1/\lambda)(-1,0)$<br>$r_3 = (0,0)$ | $(\mu_1, \mu_2, \lambda)$<br>$(\mu_1, \lambda, \mu_2)$<br>$(\lambda, \mu_1, \mu_2)$ |

4.3. *Examples and numerical results.* In this section we study two examples: the individual buffer overflow for two-node tandem Jackson network and total population overflow for $d$-node tandem Jackson network.

4.3.1. *Two-node tandem networks with individual buffer overflow.* Consider the two-node tandem queueing ($d=2$) networks with $\Theta = (\lambda, \mu_1, \mu_2)$, and the quantity of interest is

$p_n \doteq \{$size of queue 1 exceeds $B_1 n$ or size of queue 2 exceeds $B_2 n$

before the system returns to empty state, starting from $0\}$.

One can think of $B_i n$ as the individual buffer size for node $i$. In the notation we just introduced, it amounts to $\Gamma = \{x \in \mathbb{R}^2_+ : x_1 \geq B_1 \text{ or } x_2 \geq B_2\}$. Assuming $\lambda + \mu_1 + \mu_2 = 1$ and $\lambda < \min\{\mu_1, \mu_2\}$, we have (following a similar argument in [10])

$$\gamma \doteq \lim_n -\frac{1}{n} \log p_n = \min_{i=1,2} B_i \log \frac{\mu_i}{\lambda}.$$

Consider piecewise affine subsolutions that take the form $\bar{W}^n \doteq \bar{W}^n_1 \wedge \bar{W}^n_2 \wedge \bar{W}^n_3$ where

$$\bar{W}^n_k(x) \doteq \langle r_k, x \rangle + 2\gamma - k\delta_n,$$

for some small positive constants $\delta_n$. The choice of $\{r_k\}$ and its corresponding change of measure $\{\bar{\Theta}^*(r_k)\}$ are given in Table 3. See also Figures 5 and 6.

It is not difficult to check that the conditions of Theorem 4.2 are satisfied with

$$\bar{\varepsilon}_n \doteq 2\log[(\mu_1 \vee \mu_2)/\lambda] \exp\{-\delta_n/\varepsilon_n\}.$$

It follows that the corresponding importance sampling estimator is asymptotically optimal if $\delta_n \to 0, \varepsilon_n/\delta_n \to 0$ and $n\varepsilon_n \to \infty$.

4.3.2. *d-node tandem networks with total population overflow.* Consider the total population overflow for a $d$-node tandem Jackson network with $d \geq 2$, that is, $\Gamma = \{x \in \mathbb{R}_d^+ : x_1 + x_2 + \cdots + x_d \geq 1\}$ and

$$p_n \doteq \mathbb{P}\{\text{network total population reaches } n \text{ before returning to } 0,$$

$$\text{starting from } 0\}.$$

Specializing to the case $d=2$ (and assuming $\mu_1 \geq \mu_2$), the results stated in this section coincide with those of Section 3. Let $\bar{\mu} \doteq \mu_1 \wedge \mu_2 \wedge \cdots \wedge \mu_d$. Assuming $\lambda < \bar{\mu}$ and $\lambda + \mu_1 + \cdots + \mu_d = 1$, we have [10]

$$\gamma \doteq \lim_n -\frac{1}{n} \log p_n = \log \frac{\bar{\mu}}{\lambda}.$$

For any fixed $n$, we consider piecewise affine subsolutions of form $\bar{W}^n = \bar{W}_1^n \wedge \cdots \wedge \bar{W}_{d+1}^n$ where

$$\bar{W}_k^n(x) \doteq \langle r_k, x \rangle + 2\gamma - k\delta_n$$

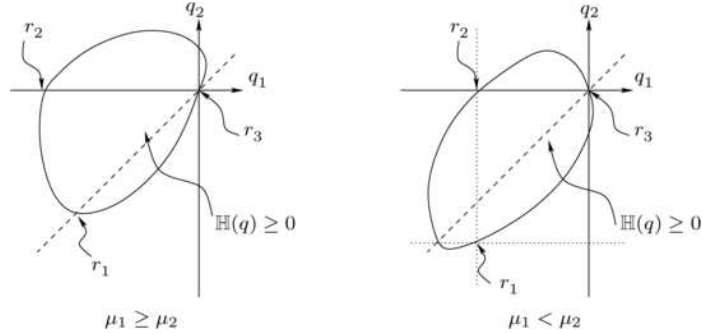

FIG. 5. *The choice of* $\{r_k\}$.

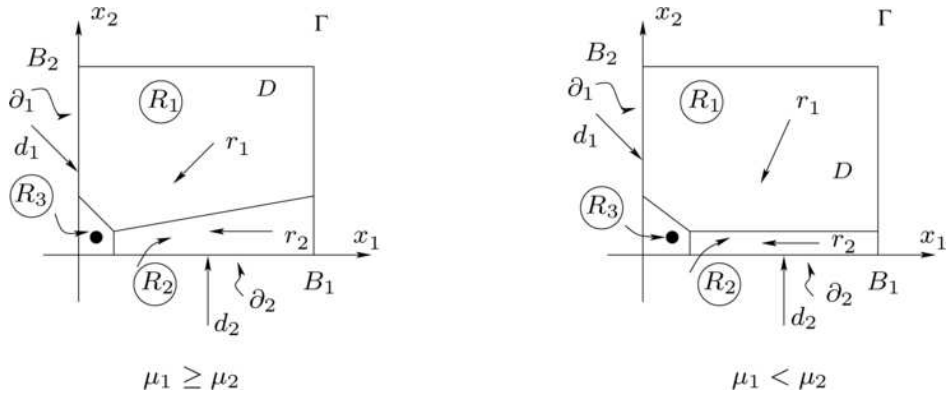

FIG. 6. *Piecewise affine function.*

for some small positive constant $\delta_n$ and

$$[r_k]_i \doteq \begin{cases} -2\gamma, & \text{if } 1 \leq i \leq d+1-k, \\ 0, & \text{otherwise,} \end{cases}$$

for $1 \leq k \leq d$ and $r_{d+1} = 0$. The change of measure corresponding to $r_k$ is

$$\bar{\Theta}^*(r_k) = \left[1 - (\mu_{d+1-k} - \bar{\mu})\frac{\bar{\mu} - \lambda}{\bar{\mu}}\right]^{-1}$$

$$\times \left(\bar{\mu}, \mu_1, \ldots, \mu_{d-k}, \frac{\lambda \mu_{d+1-k}}{\bar{\mu}}, \mu_{d+2-k}, \ldots, \mu_d\right)$$

for $1 \leq k \leq d$, and

$$\bar{\Theta}^*(r_{d+1}) = \Theta = (\lambda, \mu_1, \ldots, \mu_d).$$

We have the following lemma, whose proof is deferred to Appendix C.

LEMMA 4.4. *The following properties hold:*

1. $\mathbb{H}(r_k) \geq 0$ for every $k$,
2. $\bar{W}^n(x) \leq 0$ for all $x \in \Gamma$,
3. if $x \in \partial_i$ is such that $D\bar{W}^n(x)$ is well defined then $\langle D\bar{W}^n(x), d_i \rangle \geq 0$,
4. if $W^{\varepsilon_n,n}$ denotes the exponential weighting of $\bar{W}^n$ with $\varepsilon_n$ as the mollification parameter, then

$$\langle DW^{\varepsilon_n,n}(x), d_i \rangle \geq -\bar{\varepsilon}_n \doteq -2\gamma \exp\{-\delta_n/\varepsilon_n\}$$

*for every $x \in \partial_i$.*

Invoking Theorem 4.2, the importance sampling schemes corresponding to $W^{\varepsilon_n,n}$ are asymptotically optimal if $\delta_n \to 0$, $\varepsilon_n/\delta_n \to 0$ and $n\varepsilon_n \to \infty$.

4.3.3. *Numerical results.* For all the simulations in this section, we set $\delta = -\varepsilon \log \varepsilon$. The justification for this choice is based on an argument analogous to that of Remark 3.7.

Consider the example of a two-node tandem queue with individual buffer overflows as presented in Section 4.3.1. For the case of $\mu_1 \geq \mu_2$, we set $\lambda =$ 0.1, $\mu_1 = 0.5$, $\mu_2 = 0.4$, and $B_1 = 0.9$, $B_2 = 1$. Simulations are generated for $n = 20, 30, 40$ with $\varepsilon = 0.01$. Each estimate consists of 20,000 replications (see Table 4). Again, for comparison the theoretical value is obtained using an iterative algorithm. For the case of $\mu_1 < \mu_2$, we set $\lambda = 0.05$, $\mu_1 = 0.35$, $\mu_2 = 0.6$, and $B_1 = 1$, $B_2 = 0.6$. We run simulations for $n = 20, 30, 40$ with $\varepsilon = 0.1$, and each estimate consists of 20,000 replications (see Table 5).

TABLE 4
*IS based on subsolutions, two-node tandem, individual buffer overflow, $\mu_1 \geq \mu_2$*

|                    | $n = 20$                      | $n = 30$                      | $n = 40$                      |
|--------------------|-------------------------------|-------------------------------|-------------------------------|
| Theoretical value  | $4.81 \times 10^{-12}$        | $3.97 \times 10^{-18}$        | $3.47 \times 10^{-24}$        |
| Estimate           | $4.83 \times 10^{-12}$        | $4.04 \times 10^{-18}$        | $3.64 \times 10^{-24}$        |
| Std. Err.          | $0.20 \times 10^{-12}$        | $0.15 \times 10^{-18}$        | $0.18 \times 10^{-24}$        |
| 95% C.I.           | $[4.43, 5.23] \times 10^{-12}$ | $[3.74, 4.34] \times 10^{-18}$ | $[3.28, 4.00] \times 10^{-24}$ |

TABLE 5
*IS based on subsolutions, two-node tandem, individual buffer overflow, $\mu_1 < \mu_2$*

|  | $n = 20$ | $n = 30$ | $n = 40$ |
| --- | --- | --- | --- |
| Theoretical value | $1.44 \times 10^{-12}$ | $4.82 \times 10^{-19}$ | $1.61 \times 10^{-25}$ |
| Estimate | $1.40 \times 10^{-12}$ | $5.01 \times 10^{-19}$ | $1.85 \times 10^{-25}$ |
| Std. Err. | $0.05 \times 10^{-12}$ | $0.29 \times 10^{-19}$ | $0.21 \times 10^{-25}$ |
| 95% C.I. | $[1.30, 1.50] \times 10^{-12}$ | $[4.43, 5.59] \times 10^{-19}$ | $[1.43, 2.27] \times 10^{-25}$ |

As for the total population overflow for general $d$-node tandem networks in Section 4.3.2, we run simulations for $d = 4$ and $d = 9$.

For $d = 4$, we set $\lambda = 0.04$, $\mu_1 = \cdots = \mu_4 = 0.24$, and run simulations for $n = 20, 25, 30$ with $\varepsilon = 0.1$. Again, each estimate consists of 20,000 replications, and the theoretical value is obtained using an iterative algorithm (see Table 6).

For $d = 9$, we set $\lambda = 0.01$, $\mu_1 = \cdots = \mu_9 = 0.11$, and run simulations for $n = 20, 25, 30$ with $\varepsilon = 0.12$. Each estimate consists of 100,000 replications (see Table 7). In this case, a benchmark value is obtained using the same dynamic importance sampling algorithm but with 10 million replications (the iterative algorithm for computing the theoretical value in the case of $d = 4$ does not work here because the state space is too large).

**5. Remarks on general queueing networks.** The subsolution approach to importance sampling can be extended to general Jackson networks and networks with more general (e.g., Markov modulated) arrival/service processes. For example, a theoretical result analogous to Theorem 4.2 that applies to general open Jackson networks appears in [15]. Such extensions, even though routine to some degree, have a few distinctions. One is that Neumann-type boundary conditions, which were adequate for tandem networks, are not sufficient anymore in general, and the more elaborate boundary Hamiltonians have to be considered instead. Another distinction is that the geometric properties of the interior and boundary Hamiltonians are much less transparent. For instance, the Markov modulated case requires solving an eigenvalue problem to obtain the Hamiltonian. Consequently, explicit formulas for the gradients of the needed affine pieces are no longer available, and must be computed numerically [9].

TABLE 6
*IS based on subsolutions, four-node tandem, total population overflow*

|  | $n = 20$ | $n = 25$ | $n = 30$ |
| --- | --- | --- | --- |
| Theoretical value | $2.04 \times 10^{-12}$ | $5.02 \times 10^{-16}$ | $1.10 \times 10^{-19}$ |
| Estimate | $2.05 \times 10^{-12}$ | $5.07 \times 10^{-16}$ | $1.08 \times 10^{-19}$ |
| Std. Err. | $0.04 \times 10^{-12}$ | $0.07 \times 10^{-16}$ | $0.03 \times 10^{-19}$ |
| 95% C.I. | $[1.97, 2.13] \times 10^{-12}$ | $[4.93, 5.21] \times 10^{-16}$ | $[1.02, 1.14] \times 10^{-19}$ |

TABLE 7
*IS based on subsolutions, nine-node tandem, total population overflow*

|  | $n = 20$ | $n = 25$ | $n = 30$ |
|---|---|---|---|
| Benchmark value | $3.18 \times 10^{-14}$ | $9.41 \times 10^{-19}$ | $2.16 \times 10^{-23}$ |
| Estimate | $2.93 \times 10^{-14}$ | $10.80 \times 10^{-19}$ | $1.98 \times 10^{-23}$ |
| Std. Err. | $0.23 \times 10^{-14}$ | $1.30 \times 10^{-19}$ | $0.30 \times 10^{-23}$ |
| 95% C.I. | $[2.47, 3.39] \times 10^{-14}$ | $[8.20, 13.10] \times 10^{-19}$ | $[1.38, 2.58] \times 10^{-23}$ |

In order to illustrate some of the ideas of these generalizations, we consider the following two-node Jackson network with feedback. Again assume Poisson arrivals with rate $\lambda$ and consecutive exponentially services with rate $\mu_i$ at node $i$. However, after being served at node 2, a job has probability $\beta$ to be returned to node 1 (see Figure 7). Note that the full two node model that includes self-feedbacks and multiple arrival streams can be treated in a completely analogous fashion, albeit with more involved computations.

Suppose that the quantity of interest is the probability of total population overflow,

$$p_n \doteq \mathbb{P}\{\text{network total population reaches } n \text{ before returning to 0,}$$
$$\text{starting from 0}\}.$$

Let $\bar{\mu} \doteq \mu_1 \wedge \mu_2$. Assuming the stability condition $\lambda < \bar{\mu}(1-\beta)$, and without loss of generality, $\lambda + \mu_1 + \mu_2 = 1$, we have [10]

$$\gamma \doteq \lim_n -\frac{1}{n} \log p_n = \log \frac{(1-\beta)\bar{\mu}}{\lambda}.$$

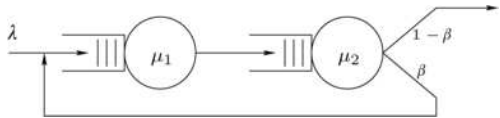

FIG. 7. *Two-node network with feedback.*

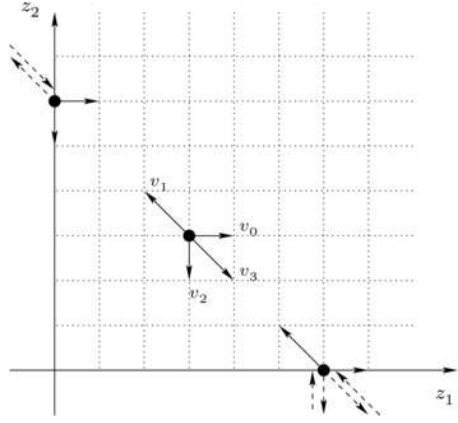

Fig. 8. *State dynamics.*

5.1. *System dynamics.* Let $Z = \{Z(k)\}$ be the embedded discrete time Markov chain that represents the queue lengths at the transition epochs of the network. Then the dynamics of $Z$ can be modeled by

$$Z(k+1) = Z(k) + \pi[Z(k), Y(k+1)],$$

where $\{Y(k)\}$ are i.i.d. random variables taking values in

$$\mathbb{V} \doteq \{v_0 = (1,0), v_1 = (-1,1), v_2 = (0,-1), v_3 = (1,-1)\},$$

and the mapping $\pi$ is defined as

$$\pi[z, y] \doteq \begin{cases} 0, & \text{if } z_1 = 0 \text{ and } y = v_1, \\ 0, & \text{if } z_2 = 0 \text{ and } y = v_2 \text{ or } v_3, \\ y, & \text{otherwise.} \end{cases}$$

Under the original probability measure $\mathbb{P}$, the distribution of $Y(k)$ is just

$$\Theta \doteq (\lambda, \mu_1, (1-\beta)\mu_2, \beta\mu_2) \in \mathcal{P}^+(\mathbb{V}).$$

See Figure 8 for an illustration of the boundary dynamics.

5.2. *The Isaacs equation and boundary Hamiltonian.* Following the argument in Section 3.5, one can write down the Isaacs equation $\mathbb{H}(DW(x)) = 0$ for $x \in D$, where

$$\mathbb{H}(p) = \sup_{\bar\Theta \in \mathcal{P}^+(\mathbb{V})} \inf_{\theta \in \mathcal{P}^+(\mathbb{V})} \left[ \langle p, \mathbb{F}(\theta) \rangle + \sum_{i=0}^{3} \theta[v_i] \log \frac{\bar\Theta[v_i]}{\Theta[v_i]} + R(\theta \| \Theta) \right]$$

with

$$\mathbb{F}(\theta) \doteq \sum_{i=0}^{3} \theta[v_i] \cdot v_i,$$

and the Dirichlet boundary condition $W(x) = 0$ for $x \in \partial_{\mathsf{e}}$.

However, as far as the boundaries $\partial_1$ and $\partial_2$ are concerned, the Neumann-type boundary condition $\langle DW(x), d_i \rangle = 0$ is not sufficient (more precisely, it is not sufficient for $\partial_2$, since the direction of constraint is not well defined on $\partial_2$). Instead one has to resort to a *boundary Hamiltonian* [6], and consequently, the boundary conditions become

$$\mathbb{H}_{\partial_i}(DW(x)) = 0 \qquad \text{for } x \in \partial_i, \ i = 1, 2,$$

where the boundary Hamiltonian $\mathbb{H}_{\partial_i}$ is defined exactly as $\mathbb{H}$ except $\mathbb{F}(\theta)$ is replaced by $\mathbb{F}_i(\theta)$ with

$$\mathbb{F}_1(\theta) = \sum_{i \neq 1} \theta[v_i] \cdot v_i, \qquad \mathbb{F}_2(\theta) = \sum_{i \neq 2,3} \theta[v_i] \cdot v_i.$$

REMARK 5.1. Proposition 4.1 can be easily applied to the interior Hamiltonian $\mathbb{H}$ and the boundary Hamiltonian $\mathbb{H}_{\partial_i}$ to show the existence of saddle points and the concavity of these Hamiltonians. The formulae for the saddle points are as follows. Let $(\bar{\Theta}^*(\cdot), \theta^*(\cdot))$ be the saddle point for $\mathbb{H}$, and $(\bar{\Theta}^*_{\partial_i}(\cdot), \theta^*_{\partial_i}(\cdot))$ be the saddle point for $\mathbb{H}_{\partial_i}$. Then

$$\bar{\Theta}^*(p) = \theta^*(p) = N(p) \cdot (\lambda e^{-p_1/2}, \mu_1 e^{(p_1-p_2)/2},$$
$$(1-\beta)\mu_2 e^{p_2/2}, \beta \mu_2 e^{(p_2-p_1)/2}),$$
$$\bar{\Theta}^*_{\partial_1}(p) = \theta^*_{\partial_1}(p) = N_1(p) \cdot (\lambda e^{-p_1/2}, \mu_1, (1-\beta)\mu_2 e^{p_2/2}, \beta \mu_2 e^{(p_2-p_1)/2}),$$
$$\bar{\Theta}^*_{\partial_2}(p) = \theta^*_{\partial_2}(p) = N_2(p) \cdot (\lambda e^{-p_1/2}, \mu_1 e^{(p_1-p_2)/2}, (1-\beta)\mu_2, \beta \mu_2),$$

where $N(p), N_i(p)$ are normalizing constants so that all these vectors are probability vectors [i.e., elements in $\mathcal{P}^+(\mathbb{V})$]. Moreover, $\mathbb{H}(p) = 2 \log N(p)$ and $\mathbb{H}_{\partial_i}(p) = 2 \log N_i(p)$.

5.3. *Piecewise affine subsolutions and mollification.* The definition of a classical subsolution is the same as Definition 3.3, except that Neumann boundary inequality $\langle DW(x), d_i \rangle \geq 0$ is replaced by $\mathbb{H}_{\partial_i}(DW(x)) \geq 0$ for $x \in \partial_i$, $i = 1, 2$.

The construction of a piecewise affine subsolution is very similar to that in Section 3.8.1 (see Figures 9 and 10). Define

$$r_1 \doteq 2\gamma(-1,-1), \qquad r_2 \doteq 2\gamma(-1,0) + 2(\gamma - a)(0,-1), \qquad r_3 \doteq (0,0),$$

where $a \in (0, \gamma]$ is given by

$$a \doteq \begin{cases} \log[\mu_1/(\mu_1 + \lambda - (1-\beta)\mu_2)], & \text{if } \mu_1 \geq \mu_2, \\ \log[\mu_1/(\lambda + \beta\mu_1)], & \text{if } \mu_1 < \mu_2. \end{cases}$$

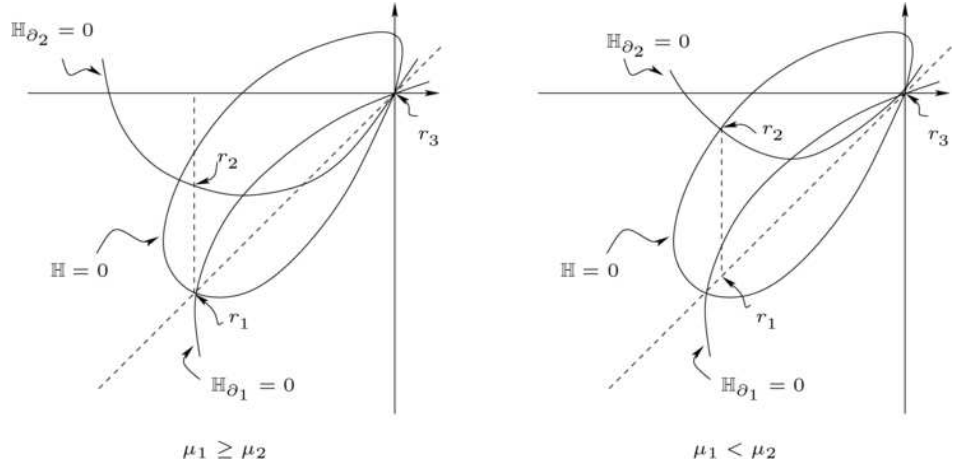

Fig. 9. *The Hamiltonians and the choice of $\{r_k\}$.*

Let $\bar{W}^\delta = \bar{W}_1^\delta \wedge \bar{W}_2^\delta \wedge \bar{W}_3^\delta$ where

$$\bar{W}_1^\delta(x) \doteq \langle r_1, x \rangle + 2\gamma - \delta,$$
$$\bar{W}_2^\delta(x) \doteq \langle r_2, x \rangle + 2\gamma - 2\delta,$$
$$\bar{W}_3^\delta(x) \doteq \langle r_3, x \rangle + 2\gamma - (1 + 2\gamma/a)\delta.$$

The exponential weighting of $\bar{W}^\delta$ with parameter $\varepsilon$ yields a smooth function

$$W^{\varepsilon,\delta}(x) \doteq -\varepsilon \log \sum_{k=1}^{3} \exp\left\{-\frac{1}{\varepsilon}\bar{W}_k^\delta(x)\right\}$$

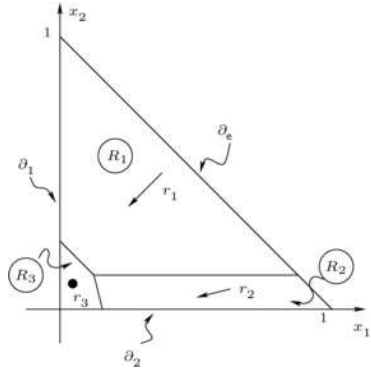

Fig. 10. *The piecewise affine subsolution.*

that satisfies

$$DW^{\varepsilon,\delta}(x) = \sum_{k=1}^{3} \rho_k^{\varepsilon,\delta}(x) r_k, \qquad \rho_i^{\varepsilon,\delta}(x) \doteq \frac{\exp\{-\bar{W}_i^{\delta}(x)/\varepsilon\}}{\sum_{k=1}^{3} \exp\{-\bar{W}_k^{\delta}(x)/\varepsilon\}}.$$

We have the following result, whose proof is deferred to Appendix C.

LEMMA 5.2. *For each $k$ we have $\mathbb{H}(r_k) \geq 0$, and the function $W^{\varepsilon,\delta}$ satisfies:*

1. $\mathbb{H}(DW^{\varepsilon,\delta}(x)) \geq 0$ *for* $x \in D$,
2. $W^{\varepsilon,\delta}(x) \leq 0$ *for* $x \in \partial_{\mathsf{e}}$,
3. *for each* $i = 1, 2$, *and* $x \in \partial_i$,

$$\mathbb{H}_{\partial_i}(DW^{\varepsilon,\delta}(x)) \geq \sum_{k=1}^{3} \rho_k^{\varepsilon,\delta}(x) \mathbb{H}_{\partial_i}(r_k) \geq -\bar{C} \exp\{-\delta/\varepsilon\}$$

*for some constant $\bar{C}$ that only depends on the system parameter $\Theta$.*

5.4. *The importance sampling scheme and its asymptotics.* Define the scaled state process $X^n(k) \doteq Z(k)/n$. Dynamic importance sampling schemes are characterized by stochastic kernels $\bar{\Theta}^n[\cdot|\cdot]$ on $\mathbb{V}$ such that $Y(k+1)$, conditional on $\{Y(1), \ldots, Y(k)\}$, has distribution $\bar{\Theta}^n[\cdot|X^n(k)] \in \mathcal{P}^+(V)$.

The importance sampling scheme based on $\bar{W}^{\varepsilon,\delta}$ is as follows. Define the stochastic kernel $\bar{\Theta}^{\varepsilon,\delta}[\cdot|\cdot]$ on $\mathbb{V}$ by

$$\bar{\Theta}^{\varepsilon,\delta}[\cdot|x] \doteq \sum_{k=1}^{3} \rho_k^{\varepsilon,\delta}(x) \bar{\Theta}^*(r_k) \qquad \text{if } x \in D$$

and

$$\bar{\Theta}^{\varepsilon,\delta}[\cdot|x] \doteq \sum_{k=1}^{3} \rho_k^{\varepsilon,\delta}(x) \bar{\Theta}^*_{\partial_i}(r_k) \qquad \text{if } x \in \partial_i.$$

Here the formulae for $\bar{\Theta}^*$ and $\bar{\Theta}^*_{\partial_i}$ can be found in Remark 5.1. We will allow $\varepsilon$ and $\delta$ to be $n$-dependent, denoted by $\varepsilon_n, \delta_n$, and let $\bar{\Theta}^n[\cdot|\cdot] \equiv \bar{\Theta}^{\varepsilon_n,\delta_n}[\cdot|\cdot]$.

Denote by $\hat{p}_n$ the corresponding importance sampling estimator. We have the following result, whose proof is very similar to that of Theorem 3.8. Indeed, in the proof of Theorem 3.8, the Neumann boundary condition is used to derive (implicitly) certain inequalities associated with boundary Hamiltonians. Such inequalities can now be obtained using Lemma 5.2. We omit the details.

THEOREM 5.3. *The importance sampling estimator $\hat{p}_n$ is asymptotically optimal if $\delta_n \to 0$, $\varepsilon_n/\delta_n \to 0$ and $n\varepsilon_n \to \infty$.*

TABLE 8
*IS based on subsolutions, two-node tandem with feedback, $\mu_1 \geq \mu_2$*

|  | $n = 20$ | $n = 30$ | $n = 40$ |
|---|---|---|---|
| Theoretical value | $9.60 \times 10^{-11}$ | $2.66 \times 10^{-16}$ | $7.27 \times 10^{-22}$ |
| Estimate | $9.31 \times 10^{-11}$ | $2.60 \times 10^{-16}$ | $7.33 \times 10^{-22}$ |
| Std. Err. | $0.17 \times 10^{-11}$ | $0.07 \times 10^{-16}$ | $0.33 \times 10^{-22}$ |
| 95% C.I. | $[8.97, 9.65] \times 10^{-11}$ | $[2.46, 2.74] \times 10^{-16}$ | $[6.67, 7.99] \times 10^{-22}$ |

One can also use a fixed pair of parameters $\varepsilon$ and $\delta$ for all $n$, which leads to a result similar to Theorem 3.6 and suggests a good choice may be to take $\delta_n = -\varepsilon_n \log \varepsilon_n$.

5.5. *Numerical results.* For all the simulations in this section, we set $\varepsilon = 0.02$ and $\delta = -\varepsilon \log \varepsilon$. For the case of $\mu_1 \geq \mu_2$, we choose $\lambda = 0.1, \mu_1 = 0.5, \mu_2 = 0.4$, and $\beta = 0.1$. Each estimate consists of 20,000 replications (see Table 8). The theoretical value is obtained using a numerical iterative algorithm.

For the case of $\mu_1 < \mu_2$, we choose $\lambda = 0.1$, $\mu_1 = 0.43$, $\mu_2 = 0.47$, and $\beta = 0.2$. Each estimate consists of 20,000 replications (see Table 9).

## APPENDIX A: A LARGE DEVIATION RESULT

In this appendix we prove a large deviation result that may be of some independent interest. Recall the definition of process $Z$ by (3.2):

$$Z(k+1) = Z(k) + \pi[Z(k), Y(k+1)],$$

where $\{Y(k)\}$ is a sequence of i.i.d. random variables taking values in $\mathbb{V} = \{v_0, v_1, v_2\}$ with distribution $\Theta = (\lambda, \mu_1, \mu_2)$. Define the hitting times

$$\sigma_n \doteq \inf\{k \geq 0 : Z_1(k) + Z_2(k) = n\},$$
$$\sigma_0 \doteq \inf\{k \geq 0 : Z_1(k) + Z_2(k) = 0\}.$$

We also let $\mathbb{Z}_n \doteq \{(z_1, z_2) \in \mathbb{Z}_+^2 : z_1 + z_2 \leq n\}$.

TABLE 9
*IS based on subsolutions, two-node tandem with feedback, $\mu_1 < \mu_2$*

|  | $n = 20$ | $n = 30$ | $n = 40$ |
|---|---|---|---|
| Theoretical value | $4.39 \times 10^{-10}$ | $2.13 \times 10^{-15}$ | $9.60 \times 10^{-21}$ |
| Estimate | $4.62 \times 10^{-10}$ | $1.91 \times 10^{-15}$ | $9.88 \times 10^{-21}$ |
| Std. Err. | $0.46 \times 10^{-10}$ | $0.13 \times 10^{-15}$ | $0.87 \times 10^{-21}$ |
| 95% C.I. | $[3.70, 5.54] \times 10^{-10}$ | $[1.65, 2.17] \times 10^{-15}$ | $[8.14, 11.64] \times 10^{-21}$ |

PROPOSITION A.1. *There exists a constant $c > 0$, which only depends on the system parameter $(\lambda, \mu_1, \mu_2)$, such that*

$$\limsup_n \sup_{z \in \mathbb{Z}_n} \frac{1}{n} \log E_z[e^{c(\sigma_n \wedge \sigma_0)}] < \infty.$$

*Here $E_z$ denotes expectation conditioned on $Z(0) = z$.*

The main difficulty in proving such a result is that the definition of $\sigma_0$ requires that the state process hit a single point, and that it is not sufficient to consider instead a small neighborhood of this point. The key idea to overcome this is to study a closely related one-dimensional process. Let $S(z) \doteq E_z[\sigma_0]$ for every $z \in \mathbb{Z}_+^2$. $S$ is finite, thanks to the stability assumption. Define the process

$$Q(k) \doteq \begin{cases} S(Z(k)), & \text{if } k \leq \sigma_0, \\ \sigma_0 - k, & \text{if } k > \sigma_0. \end{cases}$$

In other words, the process $Q$ is random until the process $Z$ hits the origin, after which $Q$ becomes deterministic and decreases by 1 each step. The scaled continuous-time piecewise affine interpolation process is just

$$Q_n(t) \doteq \frac{1}{n} Q(\lfloor nt \rfloor) + \frac{nt - \lfloor nt \rfloor}{n} [Q(\lfloor nt \rfloor + 1) - Q(\lfloor nt \rfloor)],$$

for $t \geq 0$.

In order to give a large deviation upper bound for the processes $\{Q_n\}$, we need the following definitions. Fix $\alpha \in \mathbb{R}$. For each $z \in \mathbb{Z}_+^2$, define

(A.1) $$h(z; \alpha) \doteq \log E_z \exp\{\alpha(Q(1) - Q(0))\},$$

(A.2) $$H(\alpha) \doteq \sup_{z \in \mathbb{Z}_+^2} h(z; \alpha).$$

Clearly, $H$ is convex since $h(z; \cdot)$ is convex for each $z$. The convex conjugate of $H$ is denoted by $L$, or,

(A.3) $$L(\beta) \doteq \sup_{\alpha \in \mathbb{R}} [\alpha \beta - H(\alpha)].$$

The function $L$ is nonnegative since $H(0) = 0$, and it will serve as a local upper rate function. For any fixed time $T \geq 0$, let $C([0,T]; \mathbb{R})$ be the Polish space of continuous functions on interval $[0, T]$ equipped with the supremum metric $\rho$. Define a mapping $I_T : C([0,T]; \mathbb{R}) \to \mathbb{R}_+ \cup \{\infty\}$ by

$$I_T(\phi) \doteq \begin{cases} \int_0^T L(\dot\phi(t))\, dt, & \text{if } \phi \text{ is absolutely continuous,} \\ \infty, & \text{otherwise,} \end{cases}$$

and denote its level set by

$$\Phi_x(s) \doteq \{\phi \in C([0,T]; \mathbb{R}) : \phi(0) = x, I_T(\phi) \leq s\}$$

for every $x \in \mathbb{R}$ and $s \geq 0$.

We have the following results, whose proofs are deferred to the end of this appendix. Proposition A.1 is a consequence of these two lemmas.

LEMMA A.2. *There exists a constant $M > 0$ such that $S(z) \leq M(z_1 + z_2)$ for every $z \in \mathbb{Z}_+^2$, and the absolute value of all increments of $\{Q(k)\}$ are uniformly bounded by $M$.*

LEMMA A.3. *Let $T > 0$ be given.*

1. *$I_T(\phi) \geq 0$ for every $\phi$, and $I_T(\phi) = 0$ if and only if $\dot{\phi}(t) \equiv -1$ for a.e. $t \in [0, T]$.*
2. *There exists a constant $K$ such that $I_T(\phi)$ is finite only if $\phi$ is Lipschitz continuous with Lipschitz constant $K$.*
3. *Given any compact set $F \subset \mathbb{R}$, the union of level sets, $\bigcup_{x \in F} \Phi_x(s)$, is compact for any $s \geq 0$. In particular, $I_T$ is lower semicontinuous.*
4. *For any $h > 0$ and $s \geq 0$, we have*

$$\limsup_n \sup_{z \in \mathbb{Z}_n} \frac{1}{n} \log \mathbb{P}_z \{\rho(Q_n, \Phi_{S(z)/n}(s)) > h\} \leq -s.$$

PROOF OF PROPOSITION A.1. Let $M$ be the constant in Lemma A.2, and $K$ be the Lipschitz constant in Lemma A.3. For any $\delta > 0$ and $T > 0$, define

$$F_T^\delta \doteq \{\phi \in C([0,T]; \mathbb{R}) : \phi(0) \in [0, M], -\delta \leq \phi \leq M + \delta,$$
$$\phi \text{ is absolutely continuous}, |\dot{\phi}| \leq K \vee M\},$$

which is a compact subset of $C([0,T]; \mathbb{R})$. It is not difficult to see that $I_T(\phi) > 0$ for any $\phi \in F_T^\delta$ if $T > M + \delta$. Indeed, suppose $I_T(\phi) = 0$. Then by Lemma A.3 we have $\phi(t) = \phi(0) - t$. If $\phi(0) \in [0, M]$ then for any $M + \delta < t \leq T$, $\phi(t) = \phi(0) - t < -\delta$. Thus $\phi \notin F_T^\delta$. It follows that, as long as $T > M + \delta$, $\min\{I_T(\phi) : \phi \in F_T^\delta\} > 0$, thanks to the lower semicontinuity of $I_T$ and the compactness of $F_T^\delta$.

Now fix an arbitrary $\delta$ (the specific value of $\delta$ is not important), and let $t_0 = M + 4\delta$. Define

$$s \doteq \tfrac{1}{2} \min\{I_{t_0}(\phi) : \phi \in F_{t_0}^{2\delta}\} > 0.$$

For any $x$ and $\phi \in \Phi_x(s)$, by Lemma A.3 again, $\phi$ is Lipschitz continuous with $|\dot{\phi}| \leq K$. However, $\Phi_x(s) \cap F_{t_0}^{2\delta} = \varnothing$ by definition of $s$. Therefore, for any $x \in [0, M]$ and $\phi \in \Phi_x(s)$, we must have

(A.4) $$\inf\{t \geq 0 : \phi(t) \notin [-2\delta, M + 2\delta]\} \leq t_0.$$

Define the stopping time

$$\tau_n^\delta \doteq \inf\{t \geq 0 : Q_n(t) \notin [-\delta, M + \delta]\}.$$

Thanks to Lemma A.2, $Q_n$ has Lipschitz continuous sample paths with $|\dot{Q}_n| \leq M$. Moreover, for any initial condition $Z(0) = z \in \mathbb{Z}_n$, Lemma A.2 implies $Q_n(0) = S(z)/n \in [0, M]$. It follows that

$$\mathbb{P}_z(\tau_n^\delta > t_0) = \mathbb{P}_z(Q_n \in F_{t_0}^\delta).$$

Thanks to equation (A.4), for every $Q_n \in F_{t_0}^\delta$, we have $\rho(Q_n, \Phi_{S(z)/n}(s)) > \delta$. Therefore,

$$\mathbb{P}_z(\tau_n^\delta > t_0) \leq \mathbb{P}_z(\rho(Q_n, \Phi_{S(z)/n}(s)) > \delta).$$

However, it follows from Lemma A.2 that $\{\sigma_n \wedge \sigma_0 > nt_0\} \subset \{\tau_n^\delta > t_0\}$ for $n \geq M/\delta$. As a consequence,

$$\limsup_n \sup_{z \in \mathbb{Z}_n} \frac{1}{n} \log \mathbb{P}_z(\sigma_n \wedge \sigma_0 > nt_0)$$
$$\leq \limsup_n \sup_{z \in \mathbb{Z}_n} \frac{1}{n} \log \mathbb{P}_z(\tau_n^\delta > t_0)$$
$$\leq \limsup_n \sup_{z \in \mathbb{Z}_n} \frac{1}{n} \log \mathbb{P}_z(\rho(Q_n, \Phi_{S(z)/n}(s)) > \delta)$$
$$\leq -s,$$

here the last inequality is by Lemma A.3. In particular,

$$\sup_{z \in \mathbb{Z}_n} \mathbb{P}_z(\sigma_n \wedge \sigma_0 \geq \lfloor nt_0 \rfloor + 1)$$
$$\leq \sup_{z \in \mathbb{Z}_n} \mathbb{P}_z(\sigma_n \wedge \sigma_0 > nt_0) \leq e^{-ns/2}$$

for $n$ big enough. Let $k_n \doteq \lfloor nt_0 \rfloor + 1$. Thanks to the Markov property, for all sufficiently large $n$ and all $j \geq 0$

$$\sup_{z \in \mathbb{Z}_n} \mathbb{P}_z(\sigma_n \wedge \sigma_0 \geq jk_n) \leq e^{-jns/2}.$$

Let $c$ be any constant such that $0 < c < s/(4t_0)$. We have, for $n$ big enough, $ck_n \leq ns/4$, which implies that

$$E_z[e^{c(\sigma_n \wedge \sigma_0)}] = \sum_{j=0}^\infty \sum_{i=jk_n}^{(j+1)k_n - 1} e^{ci} \mathbb{P}_z(\sigma_n \wedge \sigma_0 = i)$$
$$\leq e^{ck_n} \sum_{j=0}^\infty e^{cjk_n} \mathbb{P}_z(jk_n \leq \sigma_n \wedge \sigma_0 \leq (j+1)k_n - 1)$$
$$\leq e^{ck_n} \sum_{j=0}^\infty e^{-j(ns/2 - ck_n)}$$

$$\leq e^{ck_n} \sum_{j=0}^{\infty} e^{-jns/4}$$

$$= e^{ck_n} \frac{1}{1-e^{-ns/4}}.$$

Therefore,

$$\limsup_{n} \sup_{z \in \mathbb{Z}_n} \frac{1}{n} \log E_z[e^{c(\sigma_n \wedge \sigma_0)}] \leq \lim_{n} \frac{ck_n}{n} + \lim_{n} \frac{1}{n} \log \frac{1}{1-e^{-ns/4}} = ct_0.$$

This completes the proof. □

It remains to show Lemmas A.2 and A.3. We begin with the following result, whose proof is a straightforward consequence of the definition of $Q(k)$ and thus omitted.

LEMMA A.4. Let $\mathcal{F}_k \doteq \sigma(Z(0), Y(1), \ldots, Y(k))$. Then

$$E_z[Q(k+1) - Q(k)|\mathcal{F}_k] = -1$$

for every $z \in \mathbb{Z}_+^2$ and every $k \geq 0$.

The next lemma is concerned with the monotonicity of the sample path with respect to the initial conditions. To be more precise, for $\bar{z}, z \in \mathbb{Z}_+^2$, we say $\bar{z} \leq z$ if the inequality holds component wise. Also for $z \in \mathbb{Z}_+^2$, denote by $Z^z$ the sample path corresponding to initial condition $z$, that is,

$$Z^z(0) = z, \qquad Z^z(k+1) = Z^z(k) + \pi[Z^z(k), Y(k+1)].$$

LEMMA A.5. Define $g: \mathbb{Z}_+^2 \to \mathbb{Z}_+$ by $g(z) = z_1 + z_2$. Given any $\bar{z}, z \in \mathbb{Z}_+^2$ such that $\bar{z} \leq z$,

$$Z^{\bar{z}}(k) \leq Z^z(k),$$
$$g(Z^z(k)) - g(Z^{\bar{z}}(k)) \leq g(z) - g(\bar{z})$$

for every $k \geq 0$.

PROOF. We use induction. The claim is trivial for $k = 0$. Assume for now that it holds for some $k \geq 0$. Introduce the following notation:

$$\Gamma \doteq \{z \in \mathbb{Z}_+^2 : z_1 > 0, z_2 > 0\},$$
$$\Gamma_1 \doteq \{z \in \mathbb{Z}_+^2 : z_1 = 0, z_2 > 0\},$$
$$\Gamma_2 \doteq \{z \in \mathbb{Z}_+^2 : z_1 > 0, z_2 = 0\}.$$

We consider the following possible scenarios separately: (i) $Z^{\bar{z}}(k) \in \Gamma$; (ii) $Z^{\bar{z}}(k) \in \Gamma_1$; (iii) $Z^{\bar{z}}(k) \in \Gamma_2$; (iv) $Z^{\bar{z}}(k) = 0$. Since the proofs for these cases are essentially the same, we choose to only present case (ii). Assume that $Z^{\bar{z}}(k) \in \Gamma_1$. Thanks to the induction hypothesis $Z^{\bar{z}}(k) \leq Z^z(k)$, we must have $Z^z(k) \in \Gamma_1 \cup \Gamma$. If $Z^z(k) \in \Gamma_1$, or $Z^z(k) \in \Gamma$ but $Y(k+1) \neq v_1$, then $\pi[Z^{\bar{z}}(k), Y(k+1)] = \pi[Z^z(k), Y(k+1)]$ and the claim holds for $k+1$. For the case where $Z^z(k) \in \Gamma$ and $Y(k+1) = v_1$, we have $Z^{\bar{z}}(k+1) = Z^{\bar{z}}(k)$ and $Z^z(k+1) = Z^z(k) + v_1 = Z^z(k) + (-1, 1)$. Since $Z_1^z(k) > 0$ and $Z_1^{\bar{z}}(k) = 0$, it follows that $Z^{\bar{z}}(k+1) \leq Z^z(k+1)$. Also note that $g(Z^z(k+1)) = g(Z^z(k))$, $g(Z^{\bar{z}}(k+1)) = g(Z^{\bar{z}}(k))$. This completes the proof. $\square$

PROOF OF LEMMA A.2. Let $\bar{M} \doteq 2S((1,0)) + 2S((0,1))$. We would like to show that for any $z \in \mathbb{Z}_+^2$ and any $i = 0, 1, 2$,

(A.5) $$|S(z + \pi[z, v_i]) - S(z)| \leq \bar{M}.$$

We can assume that $\pi[z, v_i] = v_i$, since otherwise there is nothing to prove. First we consider the case $i = 2$. Let $\bar{z} \doteq z + v_2 = (z_1, z_2 - 1) \leq z$. Define stopping times $T^z \doteq \inf\{k \geq 0 : Z^z(k) = 0\}$ and $T^{\bar{z}} \doteq \inf\{k \geq 0 : Z^{\bar{z}}(k) = 0\}$. Thanks to Lemma A.5, we have $Z^{\bar{z}}(k) \leq Z^z(k)$ for any $k \geq 0$, which implies $T^{\bar{z}} \leq T^z$. By the same lemma, $g(Z^z(k)) - g(Z^{\bar{z}}(k)) \leq g(z) - g(\bar{z}) = 1$ for every $k$. In particular, for $k = T^{\bar{z}}$, it yields $g(Z^z(T^{\bar{z}})) \leq 1$. Therefore $Z^z(T^{\bar{z}}) \in \{(0,0), (1,0), (0,1)\}$. Now the strong Markov property yields

$$S(z) = S(\bar{z}) + \mathbb{P}\{Z^z(T^{\bar{z}}) = (1,0)\}S((1,0)) + \mathbb{P}\{Z^z(T^{\bar{z}}) = (0,1)\}S((0,1)).$$

Thus $|S(z) - S(\bar{z})| \leq S((1,0)) + S((0,1)) \leq \bar{M}/2$. The proof for the case $i = 0$ is almost verbatim. For $i = 1$, $z + v_i = z + (-1, 1)$. One can use the same argument to prove that $|S(z) - S(z + (-1, 0))| \leq \bar{M}/2$ and $|S(z + (-1, 0)) - S(z + (-1, 1))| \leq \bar{M}/2$, and then use the triangle inequality to show $|S(z + v_1) - S(z)| \leq \bar{M}$. We omit the details.

It follows from (A.5) that the increment of $\{Q(k)\}$ is uniformly bounded by $\bar{M}$ [note that $\bar{M} \geq 1$ trivially since $S(z) \geq 1$ for every $z \neq 0$]. Now for every $z \in \mathbb{Z}_+^2$, we can write $S(z)$ as

$$S(z) = [S(z) - S((0, z_1 + z_2))] + [S((0, z_1 + z_2)) - S((0,0))]$$
$$= \sum_{i=0}^{z_1-1} [S(z + iv_1) - S(z + (i+1)v_1)]$$
$$+ \sum_{i=0}^{z_1+z_2-1} [S(z + z_1 v_1 - iv_2) - S(z + z_1 v_1 - (i+1)v_2)].$$

Thanks to (A.5) again, the absolute value of each summand is bounded by $\bar{M}$. Thus $S(z) \leq \bar{M}(2z_1 + z_2)$. Taking $M \doteq 2\bar{M}$ completes the proof. $\square$

Proof of Lemma A.3. Clearly $h(z;\cdot)$ is convex for every $z$, $H$ is convex, and $H(0) = 0$. Let $M$ be the uniform bound on the increments of $\{Q(k)\}$ given by Lemma A.2. It follows that $|h(z;\alpha)| \leq M|\alpha|$. Therefore $|H(\alpha)| \leq M|\alpha|$ for every $\alpha$, whence $H$ is Lipschitz continuous (thanks to its convexity).

We claim that $H$ is differentiable at $\alpha = 0$ and $H'(0) = -1$. Indeed, since $h(z;\alpha)$ is differentiable with respect to $\alpha$ and $h(z;0) = 0$, we have

$$\frac{H(\alpha)}{\alpha} = \sup_{z \in \mathbb{Z}_+^2} \frac{h(z;\alpha)}{\alpha} = \sup_{z \in \mathbb{Z}_+^2} D_\alpha h(z;\alpha[z]),$$

where $\alpha[z]$ is some number between 0 and $\alpha$. But Lemma A.4 implies that $D_\alpha h(z;0) = E_z[Q(1) - Q(0)] = -1$, while Lemma A.2 and simple algebra yield that $|D_{\alpha\alpha} h(z;\alpha)| \leq \bar{K}$ for some constant $\bar{K}$ and for every $z$ and $\alpha$. Therefore,

$$\left|\frac{H(\alpha)}{\alpha} + 1\right| \leq \bar{K}|\alpha|.$$

Letting $\alpha \to 0$, it follows that $H$ is differentiable at $\alpha = 0$ with $H'(0) = -1$.

The convexity of $H$ and $H(0) = 0$ imply that $L$, defined by (A.3), is convex and nonnegative. The Lipschitz continuity of $H$ implies that $L$ takes value infinity outside a compact set. Lastly, $H'(0) = -1$ imply that $L(\beta) = 0$ if and only if $\beta = -1$. Parts 1 and 2 of Lemma A.3 follow from these properties of $L$. The rest of the lemma follows from Theorem 4.1 of [5] and that

$$L(\beta) \leq l(z;\beta) \doteq \sup_\alpha [\alpha\beta - h(z;\alpha)].$$

This completes the proof. □

## APPENDIX B: PROOF OF MAIN THEOREMS

We put the proofs of Theorem 3.6 and Theorem 3.8 together in this appendix. These proofs are, in essence, verification type arguments.

Lemma B.1. *The function $W^{\varepsilon,\delta}$ as defined in (3.13) satisfies the following properties:*

1. $\mathbb{H}(DW^{\varepsilon,\delta}(x)) \geq 0$ *for all* $x \in D$.
2. $W^{\varepsilon,\delta}(x) \leq 0$ *for all* $x \in \partial_{\mathsf{e}}$.
3. $\langle DW^{\varepsilon,\delta}(x), d_i \rangle \geq -2\gamma \exp\{-\delta/\varepsilon\}$ *for every* $x \in \partial_i$.
4. *There exists a constant $C$ which only depends on the system parameter $(\lambda, \mu_1, \mu_2)$, such that*

$$\left|\frac{\partial^2 W^{\varepsilon,\delta}(x)}{\partial x_i \, \partial x_j}\right| \leq \frac{C}{\varepsilon}$$

*for every $x \in \bar{D}$ and every $i,j$.*

PROOF. Thanks to (3.14), the concavity of $\mathbb{H}$ (Proposition 3.2), and that $\mathbb{H}(r_k) \geq 0$, it follows that

$$\mathbb{H}(DW^{\varepsilon,\delta}(x)) = \mathbb{H}\left(\sum_{k=1}^3 \rho_k^{\varepsilon,\delta}(x) r_k\right) \geq \sum_{k=1}^3 \rho_k^{\varepsilon,\delta}(x) \mathbb{H}(r_k) \geq 0.$$

By Lemma 3.5 we have $\bar{W}^{\varepsilon,\delta}(x) \leq \bar{W}^\delta(x)$. But $\bar{W}^\delta(x) \leq 0$ for $x \in \partial_{\mathsf{e}}$ by definition, and so the second claim follows.

Since $\langle r_1, d_1 \rangle = \langle r_3, d_1 \rangle = 0$ and $\langle r_2, d_1 \rangle = -2\gamma$, we have $\langle DW^{\varepsilon,\delta}(x), d_1 \rangle = -2\gamma \rho_2^{\varepsilon,\delta}(x)$. For $x \in \partial_1$, thanks to (3.15) and (3.12), we have

$$\rho_2^{\varepsilon,\delta}(x) \leq \frac{\exp\{-\bar{W}_2^\delta(x)/\varepsilon\}}{\exp\{-\bar{W}_3^\delta(x)/\varepsilon\}} = \exp\{-\delta/\varepsilon\}.$$

Once can treat $x \in \partial_2$ in an analogous fashion. This completes part 3.

Denote by $e_i$ the standard $i$th unit vector. It follows easily from (3.14) and (3.15) that

$$\frac{\partial^2 W^{\varepsilon,\delta}(x)}{\partial x_i \partial x_j} = \sum_{k=1}^3 \frac{\partial \rho_k^{\varepsilon,\delta}(x)}{\partial x_j} \langle r_k, e_i \rangle,$$

$$\frac{\partial \rho_k^{\varepsilon,\delta}(x)}{\partial x_j} = \frac{1}{\varepsilon} \cdot \rho_k^{\varepsilon,\delta}(x) \left[ -\langle r_k, e_j \rangle + \sum_{m=1}^3 \rho_m^{\varepsilon,\delta}(x) \langle r_m, e_j \rangle \right].$$

The last claim follows readily since $\rho_k^{\varepsilon,\delta}(x)$ is bounded between 0 and 1. $\square$

We now define a few functions that are closely related to the interior and boundary Hamiltonians. For each $\alpha \geq 0$ and $\bar{\Theta}, \theta \in \mathcal{P}^+(\mathbb{V})$, let

$$\bar{L}(\alpha, p; \bar{\Theta}, \theta) \doteq (1+\alpha)\langle p, \mathbb{F}(\theta)\rangle + (1+2\alpha)\sum_{i=0}^2 \theta[v_i] \log \frac{\bar{\Theta}[v_i]}{\Theta[v_i]} + R(\theta\|\Theta).$$

Similarly, for each $j = 1, 2$, let $\mathbb{F}_j(\theta) = \sum_{i \neq j} \theta[v_i] \cdot v_i$ and

$$\bar{L}_j(\alpha, p; \bar{\Theta}, \theta) \doteq (1+\alpha)\langle p, \mathbb{F}_j(\theta)\rangle + (1+2\alpha)\sum_{i=0}^2 \theta[v_i] \log \frac{\bar{\Theta}[v_i]}{\Theta[v_i]} + R(\theta\|\Theta).$$

LEMMA B.2. *Let $p \in \mathbb{R}^2$ such that $\mathbb{H}(p) \geq 0$. Then for any $\alpha \geq 0$,*

$$\inf_{\theta \in \mathcal{P}^+(\mathbb{V})} \bar{L}(\alpha, p; \bar{\Theta}^*(p), \theta) \geq 0,$$

*where $\bar{\Theta}^*(p)$ is as defined in Proposition 3.2.*

PROOF. By definition of $\bar{L}$, (3.7), and Proposition 3.2, it follows that

$$\bar{L}(\alpha, p; \bar{\Theta}^*(p), \theta) = \bar{L}(0, p; \bar{\Theta}^*(p), \theta) + 2\alpha \log N(p)$$
$$= \bar{L}(0, p; \bar{\Theta}^*(p), \theta) + \alpha \mathbb{H}(p).$$

However, thanks to Proposition 3.2 again, we have

$$\inf_{\theta \in \mathcal{P}^+(\mathbb{V})} \bar{L}(0, p; \bar{\Theta}^*(p), \theta) = \bar{L}(0, p; \bar{\Theta}^*(p), \theta^*(p)) = \mathbb{H}(p).$$

This completes the proof. □

PROOF OF THEOREM 3.6. To ease exposition, we adopt the notation $W = W^{\varepsilon, \delta}$, $\rho_k = \rho_k^{\varepsilon, \delta}$, and set $\bar{\varepsilon} \doteq 2\gamma \exp\{-\delta/\varepsilon\}$. Fix any $\alpha > 0$. We claim that

(B.1) $$\inf_{\theta \in \mathcal{P}^+(\mathbb{V})} \bar{L}(\alpha, DW(x); \bar{\Theta}^n[\cdot|x], \theta) \geq 0.$$

Indeed, thanks to the definition of $\bar{L}$, the concavity of the logarithmic function, and that $DW(x) = \sum \rho_k(x) r_k$, $\Theta^n[\cdot|x] = \sum \rho_k(x) \bar{\Theta}^*(r_k)$, we have

$$\bar{L}(\alpha, DW(x); \bar{\Theta}^n[\cdot|x], \theta) \geq \sum_{k=0}^{2} \rho_k(x) \bar{L}(\alpha, r_k; \bar{\Theta}^*(r_k), \theta).$$

Inequality (B.1) follows readily since $\mathbb{H}(r_k) \geq 0$ and Lemma B.2. Note that for every $x \in \partial_j \cap \bar{D}_n$, thanks to (B.1),

$$\bar{L}_j(\alpha, DW(x); \bar{\Theta}^n[\cdot|x], \theta)$$
$$= \bar{L}(\alpha, DW(x); \bar{\Theta}^n[\cdot|x], \theta) - (1+\alpha)\theta[v_j] \cdot \langle DW(x), v_j \rangle$$
$$\geq -(1+\alpha)\theta[v_j] \cdot \langle DW(x), v_j \rangle.$$

Recalling that $d_j = -v_j$, by Lemma B.1 we arrive at

(B.2) $$\inf_{\theta \in \mathcal{P}^+(\mathbb{V})} \bar{L}_j(\alpha, DW(x); \bar{\Theta}^n[\cdot|x], \theta) \geq -(1+\alpha)\bar{\varepsilon}.$$

We now show that inequalities (B.1) and (B.2) imply that for every $x \in \bar{D}_n$

(B.3) $$\inf_{\theta \in \mathcal{P}^+(\mathbb{V})} \left\{ \sum_{i=0}^{2} (1+\alpha) n \left[ W\left(x + \frac{1}{n}\pi(x, v_i)\right) - W(x) \right] \cdot \theta[v_i] \right.$$
$$+ (1+2\alpha) \sum_{i=0}^{2} \theta[v_i] \log \frac{\bar{\Theta}^n[v_i|x]}{\Theta[v_i]} + R(\theta \| \Theta) \right\}$$
$$\geq -(1+\alpha)\left[\frac{C}{n\varepsilon} + \bar{\varepsilon}\right],$$

where $C$ is a constant that only depends on the system parameter $(\lambda, \mu_1, \mu_2)$. To this end, consider separately the cases $x \in D_n$ (interior) and $x \in \partial_j \cap \bar{D}_n$ (boundary). For $x \in D_n$, $\pi(x, v_i) \equiv v_i$. Therefore, by Taylor's expansion,

$$n\left[W\left(x + \frac{1}{n}v_i\right) - W(x)\right] \cdot \theta[v_i]$$
$$= \langle DW(x), v_i \rangle \cdot \theta[v_i] + \frac{1}{2n} \langle v_i, D^2 W(\bar{x}_i) v_i \rangle \cdot \theta[v_i],$$

where $\bar{x}_i$ is some point on the line connecting $x$ and $x + v_i$. Thanks to Lemma B.1, the definition of $\mathbb{F}$ [see (3.7)], and that $\|v_i\|^2 \leq 2$, we have

$$\sum_{i=0}^{2} n\left[W\left(x + \frac{1}{n}v_i\right) - W(x)\right] \cdot \theta[v_i] \geq \langle DW(x), \mathbb{F}(\theta) \rangle - \frac{C}{n\varepsilon}.$$

This and (B.1) immediately lead to (B.3). The case of $x \in \partial_j \cap \bar{D}_n$ is similar, except now that $\pi(x, v_i) = v_i$ if $i \neq j$ and $\pi(x, v_j) = 0$. We omit the details.

Applying the relative entropy representation (Remark 3.1) to the left-hand side of (B.3) and adopting the notation $\beta_n \doteq C/(n\varepsilon) + \bar{\varepsilon}$, we have

$$e^{-(1+\alpha)\beta_n} \cdot \sum_{i=0}^{2} e^{-(1+\alpha)n[W(x+\pi(x,v_i)/n) - W(x)]} \left(\frac{\Theta[v_i]}{\bar{\Theta}^n[v_i|x]}\right)^{1+2\alpha} \cdot \Theta[v_i] \leq 1$$

for every $x \in \bar{D}_n$. Recalling the definition of $X^n$ in (3.4), this display implies that the process $M = \{M(k) : k \geq 0\}$, where

$$M(k) \doteq e^{-(1+\alpha)\beta_n k} e^{-(1+\alpha)nW(X^n(k))} \left(\prod_{j=0}^{k-1} \frac{\Theta[Y(j+1)]}{\bar{\Theta}^n[Y(j+1)|X^n(j)]}\right)^{1+2\alpha},$$

is a supermartingale under the original probability measure $\mathbb{P}$. Thanks to the optional sampling theorem and the nonnegativity of $M$,

$$E^{\mathbb{P}} M(T_n \wedge T_0) \leq E^{\mathbb{P}} M(0) = e^{-(1+\alpha)nW(0)}.$$

Since $\hat{p}_n = \hat{p}_n \cdot 1_{\{T_n < T_0\}}$ and $W(x) \leq 0$ for $x \in \partial_\mathsf{e}$,

$$M(T_n \wedge T_0) \geq M(T_n) \cdot 1_{\{T_n < T_0\}}$$
$$= e^{-(1+\alpha)\beta_n T_n} e^{-(1+\alpha)nW(X^n(T_n))} \hat{p}_n^{1+2\alpha}$$
$$\geq e^{-(1+\alpha)\beta_n T_n} \hat{p}_n^{1+2\alpha}.$$

It follows that

$$E^{\mathbb{P}}[e^{-(1+\alpha)\beta_n T_n} \hat{p}_n^{1+2\alpha}] \leq e^{-(1+\alpha)nW(0)}.$$

By Hölder's inequality,

[2nd moment of $\hat{p}_n$]
$$= E^{\mathbb{P}}[\hat{p}_n]$$
$$\leq E^{\mathbb{P}}[e^{-(1+\alpha)\beta_n T_n} \hat{p}_n^{1+2\alpha}]^{1/(1+2\alpha)} \cdot E^{\mathbb{P}}[e^{((1+\alpha)/(2\alpha))\beta_n T_n} \cdot 1_{\{T_n < T_0\}}]^{2\alpha/(1+2\alpha)}$$
$$\leq e^{-((1+\alpha)/(1+2\alpha))nW(0)} \cdot E^{\mathbb{P}}[e^{((1+\alpha)/(2\alpha))\beta_n (T_n \wedge T_0)}]^{2\alpha/(1+2\alpha)},$$

which yields

(B.4)
$$\liminf_n -\frac{1}{n}\log[\text{2nd moment of } \hat{p}_n]$$
$$\geq \frac{1+\alpha}{1+2\alpha}W(0) - \frac{2\alpha}{1+2\alpha}\limsup_n \frac{1}{n}\log E^{\mathbb{P}}[e^{((1+\alpha)/(2\alpha))\beta_n (T_n \wedge T_0)}].$$

Let $c$ be the constant in Proposition A.1, and let

$$\bar{C} \doteq \limsup_n \sup_{x \in \tilde{D}_n} \frac{1}{n} \log E^{\mathbb{P}}_x[e^{c(T_n \wedge T_0)}].$$

It follows immediately from Proposition A.1 that $\bar{C}$ is finite. Note that (B.4) holds for any $\alpha > 0$. In particular, it holds for $\alpha \doteq \bar{\varepsilon}/c$. With this choice of $\alpha$, we have

$$\frac{1+\alpha}{2\alpha}\beta_n = \frac{1+\alpha}{2\alpha}\frac{C}{n\varepsilon} + \frac{\bar{\varepsilon}}{2} + \frac{c}{2}.$$

Therefore, if $\bar{\varepsilon} < c$, then for $n$ big enough,

$$\frac{1+\alpha}{2\alpha}\beta_n < c.$$

It follows from (B.4) and $W(0) \leq 2\gamma$ that

$$\liminf_n -\frac{1}{n}\log[\text{2nd moment of } \hat{p}_n] \geq \frac{1+\alpha}{1+2\alpha}W(0) - \frac{2\alpha}{1+2\alpha}\bar{C}$$
$$= W(0) - \bar{\varepsilon}\frac{1}{c+2\bar{\varepsilon}}[W(0) + 2\bar{C}]$$
$$\geq W(0) - \bar{\varepsilon}\frac{1}{c}[2\gamma + 2\bar{C}].$$

It follows from Lemma 3.5 that

(B.5) $$W(0) = W^{\varepsilon,\delta}(0) \geq \bar{W}^\delta(0) - 3\varepsilon = 2\gamma - 3\delta - 3\varepsilon.$$

Recall that $\bar{\varepsilon} = 2\gamma \exp\{-\delta/\varepsilon\}$. We conclude the proof by setting $A = 2\gamma[2\gamma + 2\bar{C}]/c$, and to enforce $\bar{\varepsilon} < c$ (which was assumed in the proof) we set $B = 1/\log(2\gamma/c)$ if $c < 2\gamma$ and $B = \infty$ if $c \geq 2\gamma$. □

PROOF OF THEOREM 3.8. It suffices to show that
$$\liminf_n -\frac{1}{n}\log[\text{2nd moment of }\hat{p}_n] \geq 2\gamma,$$
since the other direction is automatic by Jensen's inequality (see Section 2). We use the notation $W^n = W^{\varepsilon_n,\delta_n}$, $\rho_k^n \doteq \rho_k^{\varepsilon_n,\delta_n}$, and $\bar{\varepsilon}_n = \exp\{-\delta_n/\varepsilon_n\}$. The same argument leading to inequality (B.4) gives that, for any strictly positive sequence $\{\alpha_n\}$,
$$\liminf_n -\frac{1}{n}\log[\text{2nd moment of }\hat{p}_n]$$
$$\geq \liminf_n \frac{1+\alpha_n}{1+2\alpha_n}W^n(0)$$
$$- \limsup_n \frac{2\alpha_n}{1+2\alpha_n}\frac{1}{n}\log E^{\mathbb{P}}[e^{((1+\alpha_n)/(2\alpha_n))\beta_n(T_n\wedge T_0)}],$$
where
$$\beta_n \doteq \frac{C}{n\varepsilon_n} + \bar{\varepsilon}_n.$$
In particular, we should choose $\alpha_n$ so that
$$\frac{1+\alpha_n}{2\alpha_n}\beta_n = c \quad \text{or} \quad \alpha_n = \frac{\beta_n}{2c-\beta_n}.$$
Note that $\alpha_n$ is strictly positive (at least for $n$ big enough) and $\alpha_n \to 0$ since $\beta_n \to 0$ by assumption. It follows that
$$\liminf_n -\frac{1}{n}\log[\text{2nd moment of }\hat{p}_n] \geq \liminf_n W^n(0).$$
However, by (B.5) $W^n(0) \geq 2\gamma - 3\delta_n - 3\varepsilon_n$. This completes the proof. □

## APPENDIX C: COLLECTION OF MISCELLANEOUS PROOFS

PROOF OF LEMMA 4.4. Clearly, $\mathbb{H}(r_{d+1}) = \mathbb{H}(0) = 0$. For $1 \leq k \leq d$, Proposition 4.1 implies that $\mathbb{H}(r_k) = 2\log N(r_k)$ where
$$\frac{1}{N(r_k)} = \bar{\mu} + \mu_1 + \cdots + \mu_{d-k} + \frac{\lambda\mu_{d+1-k}}{\bar{\mu}} + \mu_{d+2-k} + \cdots + \mu_d.$$
In order to show $\mathbb{H}(r_k) \geq 0$ or $N(r_k) \geq 1$, it suffices to show that
$$\bar{\mu} + \lambda\mu_{d+1-k}/\bar{\mu} \leq \lambda + \mu_{d+1-k},$$
or equivalently, $(\mu_{d+1-k} - \bar{\mu})(\bar{\mu} - \lambda) \geq 0$, which directly follows from the assumptions. Furthermore, for $x \in \Gamma$ we have
$$\bar{W}^n(x) \leq \bar{W}_1^n(x) = -2\gamma(x_1 + x_2 + \cdots + x_d) + 2\gamma - \delta \leq -\delta < 0.$$

Now assume $x \in \partial_i$ for some $1 \leq i \leq d$. Suppose $D\bar{W}^n(x)$ is well defined, or equivalently, $\bar{W}_1^n(x) \wedge \cdots \wedge \bar{W}_{d+1}^n(x) = \bar{W}_{k^*}^n(x)$ for some unique $k^*$. In this case, $D\bar{W}^n(x) = r_{k^*}$. In order to show $\langle r_{k^*}, d_i \rangle \geq 0$, note that

$$\text{(C.1)} \qquad \langle r_k, d_i \rangle = \begin{cases} -2\gamma, & \text{if } k+i = d+1, \\ 0, & \text{otherwise.} \end{cases}$$

Thus it suffices to show that $k^* \neq d+1-i$. This is true since the definitions of $\{r_k\}$ and $x_i = 0$ imply

$$\begin{aligned}
\bar{W}_{d+2-i}^n(x) &= \langle r_{d+2-i}, x \rangle + \gamma - (d+2-i)\delta_n \\
&= \langle r_{d+1-i}, x \rangle + \gamma - (d+2-i)\delta_n \\
&= \bar{W}_{d+1-i}^n(x) - \delta \\
&< \bar{W}_{d+1-i}^n(x).
\end{aligned}$$

It remains to show that $\langle DW^{\varepsilon_n,n}(x), d_i \rangle \geq -2\gamma \exp\{-\delta_n/\varepsilon_n\}$ for $x \in \partial_i$. Since $DW^{\varepsilon_n,n}(x) = \sum_{k=1}^{d+1} \rho_k^{\varepsilon_n,n}(x) r_k$ and for $x \in \partial_i$,

$$\rho_{d+1-i}^{\varepsilon_n,n}(x) \leq \frac{\exp\{-\bar{W}_{d+1-i}^n(x)/\varepsilon_n\}}{\exp\{-\bar{W}_{d+2-i}^n(x)/\varepsilon_n\}} = \exp\{-\delta_n/\varepsilon_n\},$$

the desired inequality follows from (C.1). □

PROOF OF LEMMA 5.2. We will only present the proof for the case $\mu_1 < \mu_2$, and omit the analogous proof for $\mu_1 \geq \mu_2$.

Assume $\mu_1 < \mu_2$ hereafter, and use the notation $W \equiv W^{\varepsilon,\delta}$ and $\rho_k = \rho_k^{\varepsilon,\delta}$. The formulae in Remark 5.1 yield

$$\mathbb{H}(r_1) = 2 \log N(r_1) = -2 \log \left[ (1-\beta)\mu_1 + \mu_1 + \beta\mu_2 + \frac{\lambda\mu_2}{\mu_1} \right].$$

By assumption $\lambda \leq (1-\beta)\mu_1$ and $\mu_1 < \mu_2$, it follows that

$$\left(\frac{\mu_2}{\mu_1} - 1\right)((1-\beta)\mu_1 - \lambda) \geq 0 \quad \text{or} \quad (1-\beta)\mu_1 + \frac{\lambda\mu_2}{\mu_1} \leq \lambda + (1-\beta)\mu_2.$$

Since $\lambda + \mu_1 + \mu_2 = 1$, we have $\mathbb{H}(r_1) \geq 0$. Similarly, we have $\mathbb{H}(r_2) = 0$ and $\mathbb{H}(r_3) = 0$. Thanks to the concavity of $\mathbb{H}$, $DW(x) = \sum_k \rho_k(x) r_k$, and $\sum_k \rho_k(x) = 1$, $\rho_k(x) \geq 0$, we have $\mathbb{H}(DW(x)) \geq 0$. As for $x \in \partial_\mathsf{e}$, we have $W(x) \leq \langle r_1, x \rangle + 2\gamma - \delta = -\delta \leq 0$.

It remains to show part 3. For $x \in \partial_i$, the concavity of $\mathbb{H}_{\partial_i}$ implies that

$$\mathbb{H}_{\partial_i}(DW(x)) \geq \sum_{k=1}^{3} \rho_k(x) \mathbb{H}_{\partial_i}(r_k) = \sum_{k=1}^{2} \rho_k(x) \mathbb{H}_{\partial_i}(r_k).$$

However, it is not difficult to check that $\mathbb{H}_{\partial_1}(r_1) \geq 0$ and $\mathbb{H}_{\partial_2}(r_2) = 0$. Therefore, we only need to show $\rho_2(x) \leq \exp\{-\delta/\varepsilon\}$ for $x \in \partial_1$ and $\rho_1(x) \leq \exp\{-\delta/\varepsilon\}$ for $x \in \partial_2$. For $x \in \partial_2$, we have $x_2 = 0$ and

$$\rho_1(x) \leq \frac{\exp\{-W_1^\delta(x)/\varepsilon\}}{\exp\{-W_2^\delta(x)/\varepsilon\}} = \exp\{-\delta/\varepsilon\}.$$

For $x = (0, x_2) \in \partial_1$, we consider two cases: $x_2 \leq x_2^*$ and $x_2 > x_2^*$ separately, where $x_2^* \doteq \delta/a$. For $x_2 \leq x_2^*$, we have

$$\rho_2(x) \leq \frac{\exp\{-W_2^\delta(x)/\varepsilon\}}{\exp\{-W_3^\delta(x)/\varepsilon\}} = \exp\left\{\frac{2(\gamma - a)}{\varepsilon}x_2 + \left(1 - \frac{2\gamma}{a}\right)\frac{\delta}{\varepsilon}\right\} \leq \exp\{-\delta/\varepsilon\}.$$

Similarly, for $x \geq x_2^*$, we have

$$\rho_2(x) \leq \frac{\exp\{-W_2^\delta(x)/\varepsilon\}}{\exp\{-W_1^\delta(x)/\varepsilon\}} = \exp\left\{-\frac{2a}{\varepsilon}x_2 + \frac{\delta}{\varepsilon}\right\} \leq \exp\{-\delta/\varepsilon\}.$$

This completes the proof. $\square$

P. DUPUIS
A. D. SEZER
H. WANG
DIVISION OF APPLIED MATHEMATICS
BROWN UNIVERSITY
PROVIDENCE, RHODE ISLAND 02912
USA
E-MAIL: dupuis@cfm.brown.edu
       sezer@cfm.brown.edu
       huiwang@cfm.brown.edu